\documentclass[a4paper,10pt]{article}

\usepackage[frenchb,british]{babel}
\usepackage{amsfonts}
\usepackage{amsmath}
\usepackage{amssymb}
\usepackage[utf8]{inputenc}
\usepackage{amsthm}
\usepackage{graphicx}
\usepackage{tikz}
\usepackage{tikz-cd}
\usepackage{hyperref}
\usepackage{amssymb}

\hypersetup{                    % parametrage des hyperliens
    colorlinks=true,                % colorise les liens
    breaklinks=true,                % permet les retours à la ligne pour les liens trop longs
    urlcolor= blue,                 % couleur des hyperliens
    linkcolor= blue,                % couleur des liens internes aux documents (index, figures, tableaux, equations,...)
    citecolor= blue               % couleur des liens vers les references bibliographiques
    }

\theoremstyle{definition}
\newtheorem{definition}{Definition}[section]
\newtheorem{thm}[definition]{Theorem}

\newtheorem{lem}[definition]{Lemma}
\newtheorem*{dem}{Proof}
\newtheorem{prop}[definition]{Proposition}
\newtheorem{cor}[definition]{Corollary}

\newtheorem{Expl}[definition]{Example}
\newtheorem{rk}[definition]{Remark}

\newcommand{\N}{\mathbb N}
\newcommand{\Z}{\mathbb Z}
\newcommand{\R}{\mathbb R}
\newcommand{\C}{\mathbb C}

\usepackage{geometry}
\usepackage{fancyhdr}

\usepackage{authblk}
\title{Controlled $K$-theory for groupoids \& applications to Coarse Geometry}
\date{}
\author{ Clément Dell'Aiera}
\affil{Department of Mathematics, University of Hawai‘i at Manoa, Honolulu, HI, 96822, U.S.A.
E-mail address: dellaiera@math.hawaii.edu}
\begin{document}
\maketitle
\begin{abstract}
We develop a generalization of quantitative $K$-theory, which we call controlled $K$-theory. It is powerful enough to study the $K$-theory of crossed product of $C^*$-algebras by action of \'etale groupoids and discrete quantum groups. In this article, we will use it to study groupoids crossed products. We define controlled assembly maps, which factorize the Baum-Connes assembly maps, and define the controlled Baum-Connes conjecture. We relate the controlled conjecture for groupoids to the classical conjecture, and to the coarse Baum-Connes conjecture. This allows to give applications to Coarse Geometry. In particular, we can prove that the maximal version of the controlled coarse Baum-Connes conjecture is satisfied for a coarse space which admits a fibred coarse embedding, which is a stronger version of a result of M. Finn-Sell.
\end{abstract}
\tableofcontents
\vspace{0.5 cm}
\textbf{Acknowledgments} This work was completed during a PhD under the supervision of H. Oyono-Oyono. The author would like to thank his advisor for very helpful discussions, comments and corrections. He also very grateful toward Rufus Willett for helpful remarks. 

%%%%%%%%%%%%%%
\section{Introduction}

Let $X$ be a countable discrete metric space with bounded geometry, i.e. such that, for every $R>0$, 
\[\sup_{x\in X} |B(x,R)|<\infty.\] 
Let $H$ be the separable Hilbert space and 
\[C_R[X] = \{T\in \mathcal L(H\otimes l^2(X)) \text{ s.t. } T_{xy} \in \mathfrak K(H) \text{ and } prop(T) < R \}\]
where $prop(T) = \sup\{d(x,y) \text{ s.t. } T_{xy} \neq 0\}$. Recall that the Roe algebra of $X$ is the $C^*$-algebra defined by
\[C^*(X) = \overline{\cup_{R>0} C_R[X]}\quad ,\] 
where the closure is taken under the operator norm. The coarse Baum-Connes conjecture asserts that
\[\mu_X : \varinjlim KK(C_0(P_d(X),\C) \rightarrow K(C^*(X))\]
is an isomorphism, where:
\begin{itemize}
\item[$\bullet$] $KK(A,B)$ denotes the operator $KK$-theory of G. Kasparov,
\item[$\bullet$] $K(B)$ denotes the operator $K$-theory of the $C^*$-algebra $B$,
\item[$\bullet$] $P_d(X)$ is the Rips complex of $X$.
\end{itemize}

The main application of the coarse Baum-Connes conjecture is the Novikov conjecture on the homotopy invariance of the higher signatures. More precisely, let $\Gamma$ be a finitely generated 
%countable discrete 
group endowed with any left-invariant metric. Such a metric is unique up to quasi-isometry, and let us denote by $|\Gamma|$ the coarse class of the underlying metric space. Denote by $B \Gamma$ the classifying space of $\Gamma$, by $[M]\in H_{dim(M)}(M,\mathbb Q)$ the fundamental class of $M$ and by $\mathcal L_M$ the $L$-class of $M$.

\begin{thm}[Descent principle]
If $B\Gamma$ has the homotopy type of a finite CW-complex and $\mu_{|\Gamma|}$ is an isomorphism, then the Novikov conjecture holds for $\Gamma$, i.e. for any $x\in H^*(B\Gamma,\mathbb Q)$, any map $f:M\rightarrow B  \Gamma $, the higher signature
\[ \sigma(M,f) = \langle \mathcal L_M\cup f^*(x) ,[M]\rangle\]
is homotopy invariant.
\end{thm}  

In \cite{Yu1}, G. Yu proved that if $X$ is of finite asymptotic dimension, then the coarse Baum-Connes conjecture holds for $X$. More generally, G. Yu proved the following theorem.

\begin{thm}[G. Yu \cite{Yu2}]
If $X$ has property $A$, then the coarse Baum-Connes conjecture holds for $X$.
\end{thm}

As property $A$ is implied by finite asymptotic dimension, the last theorem is more powerful. Still, the two results are quite different in their proofs. Whereas property $A$ entails the existence of so called Dirac and Dual-Dirac elements in $KK$-theory, yielding the result, the proof in the setting of finite asymptotic dimension relies on a analog of Mayer-Vietoris decompositions on the Roe algebras. This proof is more geometric in nature, and more elementary. Quantitative $K$-theory was introduced in \cite{OY2} by H. Oyono-Oyono and G. Yu in order to broaden the domain of validity of this strategy.\\ 

In the author's thesis was introduced controlled $K$-theory. This slight generalization of quantitative $K$-theory relies on a new definition of filtered $C^*$-algebras, which allows to treat more examples of $C^*$-algebras with the controlled $K$-theory groups. We define controlled assembly maps in the setting of Roe algebras and of crossed products of $C^*$-algebras by étale groupoids. These assembly maps take values in the controlled $K$-theory, and should enjoy more stability properties than the usual assembly maps. The latter is the object of future work. \\

The controlled assembly maps induces the assembly maps in $K$-theory. We study in more details this phenomenon, which gives what we call quantiative statements (theorems \ref{Quant1},\ref{Quant2} \& \ref{UniformQS}). These theorems relate the classical Baum-Connes conjecture for an étale groupoid to its controlled analog. As a byproduct, we prove in lemma \ref{prod} an interesting result on the $K$-homology of a finite $G$-simplicial complex with values in an infinite product of stable $C^*$-algebras.\\

Following the route of \cite{SkTuYu}, we show that these controlled assembly maps are related by the coarse groupoid defined by G. Skandalis, J-L. Tu and G. Yu. More precisely, out of any coarse space $X$, one can construct an étale groupoid $G(X) \rightrightarrows \beta X$, such that the coarse assembly map $\mu_{X,B}$ is equivalent to the assembly map

\[\mu_{G(X),l^\infty_B} : \varinjlim KK(C_0(P_E(G),B) \rightarrow K(l^\infty_B\rtimes_r G(X))\]
where :
\begin{itemize}
\item[$\bullet$] $\beta X$ is the Stone-\v{C}ech compactification of $X$,
\item[$\bullet$] the inductive limit is taken over the compact subsets $E\subseteq G$,
\item[$\bullet$] $l^\infty_B$ is the $G(X)$-algebra $l^\infty (X,B\otimes \mathfrak K)$ and $l^\infty_B\rtimes_r G(X)$ is the associated reduced crossed product,
\item[$\bullet$] $P_E(G)$ is the Rips complex of $G$.
\end{itemize}
These results are implied by their analog in controlled $K$-theory that we prove (theorem \ref{BCCeq}). As a corollary, we prove that any coarse space which admits a fibred coarse embedding into Hilbert space satisfies the controlled maximal coarse Baum-Connes conjecture (theorem \ref{fibred}). This is a stronger version of a result of M. Finn-Sell \cite{FinnSellFibred}. Recall that the notion of fibred embedding into Hilbert space is weaker than embedding into Hilbert space. For instance, some box spaces of $SL(2,\Z)$ are expanders, hence cannot embed into Hilbert space, but admits a fibred embedding.\\ 

The article follows the following plan. The second section present an overview of controlled $K$-theory in the setting of $C^*$-algebras filtered by what we call a coarse structure. In the third and fourth sections, we build assembly maps with values in these controlled $K$-groups, which factorizes the usual assembly maps, in the case of étale groupoids and coarse spaces respectively. The last section is devoted to applications of these results in Coarse Geometry. 

\subsection{Preliminaries on groupoids}

Recall the following definition.

\begin{definition}
An étale groupoid is given by two topological spaces, the space of arrows $G$ and the space of units $G^{(0)}$ endowed with:
\begin{itemize}
\item[$\bullet$] continuous maps $s,r : G \rightrightarrows G^{(0)}$ which are local homeomorphisms,
\item[$\bullet$] a topological embedding $e: G^{(0)}\rightarrow G$ called the unit map, and a continuous involution $inv : G\rightarrow G; g\mapsto g^{-1}$ called the inverse map,
\item[$\bullet$] a multiplication map $G\times_{s,r}G\rightarrow G; (g,g')\mapsto gg'$ such that $(gg')g'' = g(g'g'')$, $gg^{-1}= e_{r(g)}$, $g^{-1}g= e_{s(g)}$
\end{itemize}
\end{definition}

\subsubsection{Actions of groupoids}

\begin{definition}
A right action of $G$ on a topological space $Z$ is given by a continuous map $p : Z \rightarrow G^{(0)}$, called the anchor map or the moment map, and a map $\alpha : Z\times_{p,r} G \rightarrow Z $ such that :
\begin{itemize}
\item[$\bullet$] $\alpha(\alpha(z,g),g') = \alpha(z, gg')$ whenever $(g,g')\in G^{(2)}$ and $p(z)=r(g)$,
\item[$\bullet$] $p(\alpha(z,g))= s(g)$
\item[$\bullet$] $\alpha(z,e_{p(z)})=z$
\end{itemize} We will use the notation $\alpha(g,z) = z.g$ when the action is clear from the context.\\ 
\end{definition}

Let $Z$ be a right $G$-space. Define :
\begin{itemize}
\item[$\bullet$] $(Z\rtimes G)^{(0)} = Z$, $Z\rtimes G = Z\times_{r,p} G$ as a topological space,
\item[$\bullet$] $u_{z}= (z,e_{p(z)})$, $ s(z,g) = z$ and $r(z,g)=z.g$, 
\item[$\bullet$] if $y = x.g$, $(x,g)(y,g')= (x,gg') $ and $(x,g)^{-1} = (x.g,g^{-1})$.
\end{itemize}

These maps define a structure of toplogical groupoid on $Z\rtimes G$. It is called the crossed product groupoid of $Z$ by $G$. It is étale if $G$ is.\\

We present now an important class of $G$-spaces called $G$-simplicial complexes. The reader is referred to \cite{TuBC2} for details. 

\begin{definition}
A map between two topological spaces $f : X\rightarrow Y$ is said to be locally injective if there exists an open cover $\mathcal U$ of $X$ such that, for all $U\in \mathcal U$, $f_{|U}$ is injective.
\end{definition}

\begin{definition} \label{Gcomplex}
Let $n\in\N$. A $G$-simplicial complex of dimension $\leq n$ is a pair $(X,\Delta)$ where :
\begin{itemize}
\item[$\bullet$] $X$ is a locally compact proper $G$-space, called the space of vertices, such that the anchor map $p : X\rightarrow G^{(0)}$ is locally injective;
\item[$\bullet$] $\Delta$ is a closed $G$-invariant subset of the space of measures on $X$, denoted $M_X$, endowed with the weak $*$-topology. Moreover, $\Delta$ contains only probability measures and satisfies :
\begin{itemize}
\item[$\bullet$] for all $\eta\in\Delta$, there exists $x\in G^ {(0)}$ such that $\text{supp }\eta \subseteq p^{-1}(x)$ and $|\text{supp }\eta|\leq n+1$,
\item[$\bullet$] if $\eta' \in \Delta$ and $\eta\in M_X$ such that $\text{supp }\eta \subseteq \text{supp }\eta'$, then $\eta\in \Delta$.
\end{itemize}
For $\eta\in \Delta$, $\text{supp }\eta$ is called a simplex, or a $j$-simplex when $|\text{supp }\eta | = j$.
\end{itemize}
The complex is said to be typed if there exists a finite space $T$ and a $G$-invariant continuous map $\tau : X\rightarrow T$ such that, for every simplex $S$, $\tau_{|S}$ is injective.  
\end{definition}

\subsubsection{Crossed products for groupoids}

We give a short review of the crossed product construction for étale groupoids. If $p : Y_0\rightarrow X$ and $q : Y_1\rightarrow X$ are two fibrations, we denote by $Y_0\times_{p,q} Y_1 = \{(y,y')\in Y_0\times Y_1 \text{ s.t. } p(y)=q(y')\}$ their fibred product. The non-commutative anologue of fibration are given by $C(X)$-algebras. 

\begin{definition}
A $C(X)$-algebra is given by a $C^*$-algebra $A$ and a non-degenerate $*$-homomorphism $\theta : C_0(X) \rightarrow Z(\mathcal M(A))$, where $\mathcal M(A)$ is the $C^*$-algebra of multipliers of $A$.
\end{definition}

If $(A,\theta)$ is a $C(X)$-algebra and $x\in X$, then the fiber over $X$ is defined as $A_x : = A/ \text{ker }(ev_x)A$. Two $C(X)$-algebras have a balanced tensor product $A\otimes_{C(X)} A'$, such that $(A\otimes_{C(X)} A')_x \cong A_x\otimes A'_x$. Notice that in the case of fibration, we get $C_0(Y_0\times_{p,q}Y_1) \cong C_0(Y_0)\times_{C(X)}C_0(Y_1)$. If $f:X_0 \rightarrow X_1$ is a continuous map, and $A$ is a $C(X_1)$-algebra, $f^* A :=C(X_0)\otimes_{C(X_1)} A$ is naturally a $C(X_0)$-algebra, called the pull-back of $A$ along $f$. For details in $C(X)$-algebras, the reader can consult \cite{blanchard} for instance. 
\\

If $A$ is a $C(G^{(0)})$-algebra, an action of the groupoid $G$ is a isomorphism of $C(G)$-algebras $\alpha : s^* A \rightarrow r^* A$ such that $\alpha_g \circ \alpha_{g'} = \alpha_{gg'}$ for every $(g,g')\in G\times_{s,r} G$. Then, define the space of continuous sections with compact support 
\[C_c(G,A) = \bigcup_{U} C_0(U)\otimes_s A\]
where $U$ runs along all open relatively compact subsets of $G$. Here $C_0(U)\otimes_s A$ denotes $C_0(U)\otimes_{C(G^{(0)})} A$, where $s$ in subscript implies that the $C(G^{(0)})$-algebra structure on $C_0(U)$ is given by $s$.\\

Endowing compact sections with convolution
\[(f_0\ast f_1)(g) = \sum_{h\in G^{r(g)}} f_0(h) \alpha_h(f_1(h^{-1}g)).\]
and involution $\overline f(g)=\alpha_g(f(g^{-1})^*)$ for every $f_0,f_1\in C_c(G,A)$, we get a $*$-algebra. 

The $A$-Hilbert module $L^2(G,A)$ is the completion of $C_c(G,A)$ under the scalar product 
\[\langle \xi ,\eta \rangle_x  = \sum_{g\in G^x} \xi(g)\overline \eta(g) \quad x\in G^{(0)} \]
and $C_c(G,A)$ is represented on $L^2(G,A)$ by $\lambda(f) \xi = f\ast \xi$, for every $ f\in C_c(G,A)$ and $\xi\in L^2(G,A)$.\\

\begin{definition}
The reduced crossed product $A\rtimes_r G$ is the $C^*$-algebra obtained by completion of $C_c(G,A)$ under the norm $||f||_r=||\lambda(f)||$.
\end{definition}

In the following section is presented the setting of controlled $K$-theory. The main example is that the set $\mathcal E$ of compact subsets of $G$ defines a coarse structure, and that $A\rtimes_r G$ is $\mathcal E$-filtered, thus allowing to apply the controlled machinery to crossed product of groupoids. \\

\subsubsection{Equivariant $KK$-theory}

We will use intensively a version of the bivariant Kasparov theory developed by P-Y. Le Gall in his thesis \cite{LeGall}, which is an equivariant $KK$-theory in the setting of groupoids, denoted by $KK^G$. Recall the following lemma.

\begin{lem}[Stabilization lemma]
Let $A$ be a $C^*$-algebra and $E$ a countably generated $A$-Hilbert module. Then, there exists an isomorphism of $A$-Hilbert modules $H_A \cong E\oplus H_A$, hence there exists a projection $p\in \mathcal L_A(H_A)$ such that $E\cong p H_A$. 
\end{lem}

We will need the following lemma, in order to put some $K$-cycles in some standard form.

\begin{rk}\label{isometry}
When we look at a $*$-homomorphism $\phi : A\rightarrow B$, we always have an isometry $V\in \mathcal L_B ( H_A\otimes_\phi B , H_B)$ defined on simple tensors as $(x_j)_j\otimes b \mapsto (\phi(x_j)b)_j$. Indeed, this map extends linearly to $H_A \odot B$, and if $x = (x_j)$ and $x'=(x'_j)$ are in $H_A$ : 
\[\langle V (x\otimes b) , V(x'\otimes b')\rangle = b^* \sum_j \phi(x_j)^* \phi(x'_j) \  b' = b^*\phi(\langle x, x' \rangle)b' = \langle x\otimes b , x'\otimes b' \rangle . \] % Determiner V^* 
This isometry can be used to explicitly describe the projection and the isomorphism appearing in the stabilization theorem in this particular example. Indeed $p = VV^*\in\mathcal L_B(H_B)$ is a projection such that $p H_B \cong H_A\otimes_\phi B $.\\

Moreover, for $T\in \mathcal L_B(H_A\otimes_\phi B)$, $Ad_V(T) = VTV^*$ defines a $*$-homomorphism $Ad_V : \mathcal L_B(H_A\otimes_\phi B)\rightarrow \mathcal L_B(H_B)$ such that $Ad_V(\mathfrak K_B(H_A\otimes_\phi B))\subseteq \mathfrak K_B(H_B)$. Indeed, notice that 
\[V\theta_{\xi,\xi'}V^* = \theta_{V\xi,V\xi'}\]
for every $\xi,\xi'\in H_A\otimes_\phi B$.\\

Composing with $\phi_*$, we get a $*$-homomorphism $\mathcal L_A(H_A)\rightarrow \mathcal L_B(H_B); T\mapsto V(T\otimes_\phi 1)V^*$ respecting compact operators in a natural way : if $\theta = (a_{ij})\in A\otimes\mathfrak K$, then $V(T\otimes_\phi 1)V^* = (\phi(a_{ij}))\in B\otimes\mathfrak K$.  
\end{rk}

\begin{lem}\label{isometryKK}
Let $\phi : B\rightarrow B'$ be a $G$-equivariant homomorphism, and $z=[H_B,\pi, T]\in KK^G(A,B)$. Let $V\in\mathcal L_{B'}(H_B\otimes B', H_{B'})$ be the isometry of the remark \ref{isometry} and $p = VV^*\in\mathcal L_{B'}(H_{B'})$. Define $\pi' : A\rightarrow \mathcal L_{B'}(H_{B'})$ as $\pi'(a) = V\pi(a)V^*$ and $T'= V(T\otimes_\phi 1)V^* + 1-p \in \mathcal L_{B'}(H_B')$. Then $(H_{B'},\pi',T')$ is a $K$-cycle and 
\[g^*(z) = [H_B\otimes B',\pi\otimes_\phi 1,\phi_*(T)]=[H_B', \pi', T']\text{ in } KK^G(A,B').\]
\end{lem} 

The last property we need to recall is decomposition property $(d)$. H. Oyono-Oyono has shown in the appendice of \cite{LaffOY} that every element of $KK^G(A,B)$ can be written as the Kasparov product of at most $d$ elements, each one coming either from a $*$-homormorphism or from a $KK$-inverse of a $*$-homomorphism.\\

\begin{definition}\label{DecompositionPropertyD}
Let $d$ be a positive integer. An element $z\in KK^G(A,B)$ is said to satisfy decomposition property $(d)$ if
\begin{itemize}
\item[$\bullet$] there exist $G$-algebras $A_0$, $A_1$, ..., $A_d$ such that $A_0=A$ and $A_d=B$, 
\item[$\bullet$] there exist elements $z_j \in KK^G(A_{j},A_{j+1})$ for $j\in\{0,..,d-1\}$ such that, either $z_j$ is induced by a $G$-morphism $A_j \rightarrow A_{j+1}$, or there exists a $G$-morphism $\phi_j : A_{j+1}\rightarrow A_j$ such that $z_j \otimes_{A_{j+1}} [\phi_j] = 1_{A_j}$ and $ [\phi_j] \otimes_{A_{j}} z_j  = 1_{A_{j+1}}$,
\end{itemize}
such that $z = z_1 \otimes_{A_1}  ... \otimes_{A_{d-1}} z_{d-1} $ holds.
\end{definition}

Then, the following theorem is true for a universal constant $d$, which does not depend on the groupoid. It will be crucial to prove that the controlled Kasparov and Roe transforms, applications to be defined later, respect the Kasparov product. 

\begin{thm}[Theorem $2.2$ \cite{LaffOY}]\label{propertyD}
Let $G$ be a locally compact groupoid with Haar system. Then, there exists a universal constant $d$ such that every element $z\in KK^G(A,B)$ has decomposition property $(d)$.
\end{thm}

Finally, let us briefly recall from \cite{TuBC} the construction of the assembly map for groupoids. \\

For any compact subset $E\subseteq G$, define $P_E(G)$ to be the space of probability measures $\nu $ with support contained in one and only one fiber $G^x$ for some $x\in G^{(0)}$, and such that if $g,g'\in \text{supp }(\nu)$, then $g'g^{-1}\in E$. We endow $P_E(G)$ with the weak-$*$ topology. It is a proper and cocompact $G$-space for left translation.\\

Every element $\eta$ is a finite probability measure on a fiber $G^x$, for some $x\in G^{(0)}$, hence can be written as a finite convex combination $\eta = \sum_{g\in G^{x}}\lambda_g(\eta)\delta_g$, where $\lambda_g(\eta)\in [0,1]$ for every $g$ and $\delta_g$ is the Dirac probability measure at $g\in G^x$. Then 
\[\mathcal L_E(g,\eta) = \lambda_{e_x}^{\frac{1}{2}}(\eta)\lambda_g^{\frac{1}{2}}(\eta)\]  
 defines a projection in $C_0(P_E(G))\rtimes_r G$ with support in $E$.

\begin{definition}\label{projection}
The assembly map for $G$ with coefficients in $B$ is defined as the inductive limit of the maps $\mu_{G,B}^{E} : KK^G(C_0(P_E(G)),B)\rightarrow K(B\rtimes_r G)$ given by
\[\mu_{G,B}^{E} (z)=[\mathcal L_E]\otimes_{C_0(Z)\rtimes G} j_G(z),\]
that is $\mu_{G,B} = \varinjlim \mu_{G,B}^{E}$ (one has to check that theses maps respects the inductive systems, which they do).\\
\end{definition}

\subsection{Coarse geometry}

Let $X$ be a countable discrete metric space with bounded geometry. The set of entourages is the set of subsets $E \subseteq X\times X$ such that $\sup d_{|E}<\infty$. We recall in this paragraph the construction of the Roe algebra of $X$ with coefficients in a $C^*$-algebra $B$. For every symmetric entourage $E\subseteq X\times X$, define $C_E[X,B]$ as the following subspace of $\mathcal L_B(H\otimes l^2(X)\otimes B)$ :
\[C_E[X,B] = \{T\in \mathcal L_B(H\otimes l^2(X)\otimes B) \text{ locally compact  s.t. supp }T\subseteq E \}.\]
It is a subspace of $\mathcal L_B(H\otimes l^2(X)\otimes B)$ which satisfies $C_E[X,B].C_{E'}[X,B]\subseteq C_{E\circ E'}[X,B]$, where $E\circ E' = EE'\cup E'E$. It is easy to see that 
\[C[X,B] = \bigcup_{E} C_E[X,B]\] 
is an involutive sub-algebra of $\mathcal L_B(H\otimes l^2(X)\otimes B)$.

\begin{definition}
Let $B$ be a $C^*$-algebra. The Roe algebra of $X$ with coefficients in $B$, denoted by $C^*(X,B)$, is the completion of $C[X,B]$ under the operator norm of $\mathcal L_B(H\otimes l^2(X)\otimes B)$. 
\end{definition} 

The following property will be useful.

\begin{thm}\label{Xfunctor}
Let $X$ be a discrete metric space with bounded geometry and $\phi : A\rightarrow B$ a $*$-homomorphism. Then there exists a $*$-homomorphism $\phi_X : C^*(X,A)\rightarrow C^*(X,B)$ extending $\phi$. Moreover, $\phi\mapsto \phi_X$ respects composition of $*$-homomorphisms.
\end{thm}

\begin{dem}
Recall that any $*$-morphism $\phi : A\rightarrow B$ induces, for any $A$-Hilbert module $E$, a $*$-morphism $\phi_* : \mathcal L_A(E)\rightarrow \mathcal L_B(E\otimes_A B)$. Now take $E$ to be $l^2(X)\otimes A$. Then, according to remark \ref{isometry}, $\eta\otimes a\otimes b\mapsto \eta \otimes\phi(a) b $ extends to an isometry $V\in \mathcal L_B(E\otimes_A B,l^2(X)\otimes B)$ which respects compact operators. The composition $Ad_V\circ\phi_*$ sends a compact operator $(T_{xy})_{x,y}$ to $(\phi(T_{xy}))_{x,y}$.\\

Hence, the linear map $T \mapsto V\phi_*(T)V^*$ maps $C_R[X,A]$ into $C_R[X,B]$, and so extends to a $*$-morphism $C^*(X,A)\rightarrow C^*(X,B)$. The composition property is clear from the construction.\\
\qed
\end{dem}

\begin{rk}
Notice that $\phi\rightarrow \phi_X$ can be defined in the same way for completely positive maps : for every completely positive map $\phi: A \rightarrow B$, there exists a completely positive map $\phi_X : C^*(X,A)\rightarrow C^*(X,B)$ extending $\phi$.
\end{rk}

%%%%%%%%%%%%%%%%%%%%%%%%%%
\section{$K$-theory controlled by a coarse structure}

In this section, we define controlled $K$-theory in more generality than in \cite{OY2}. The goal is to develop controlled $K$-theory to a broader setting that what was used until now. We start by the definition of a coarse structure, which will be the index set of the filtration of filtered $C^*$-algebras in our sense. The setting of controlled $C^*$-algebras allows one to extend the "quantitative philosophy" developed by H. Oyono-Oyono and G. Yu to the realm of crossed-product of $C^*$-algebras by actions of groupoids and discrete quantum groups. This last example will not be explained in this article, the interested reader can consult the author's thesis \cite{DellAieraThesis} for details. 

\begin{definition}
A coarse structure $\mathcal E$ is a lattice which is an abelian semi-group. %such that $\forall E,E'\in \mathcal E$, $E\leq E'E$. 
Recall that a lattice is a poset for which every pair $(E,E')$ admits a supremum $E\vee E'$ and an infimum $E\wedge E'$.
\end{definition}

\begin{definition}
A $C^*$-algebra $A$ is said to be $\mathcal E$-filtered if there exists a coarse structure $\mathcal E$ and, for every $E\in \mathcal E$, linear subspaces $A_E$ of $A$ such that :
\begin{itemize}
\item[$\bullet$] if $E \leq E'$, then $A_E\subseteq A_{E'}$, and the inclusion $\phi_E^{E'}: A_E\hookrightarrow A_{E'}$ induces an inductive system of linear spaces,
\item[$\bullet$] $A_E$ is stable by involution,
\item[$\bullet$] for all $E,E'\in\mathcal E$, $A_E.A_{E'}\subseteq A_{EE'}$,
\item[$\bullet$] the union of subspaces is dense in $A$, i.e. $\overline{\cup_{E\in\mathcal E}A_E} = \varinjlim A_E = A$.
\item[$\bullet$] if $A$ is unital, we impose that $1\in A_E,\forall E\in\mathcal E$.
\end{itemize}
\end{definition}

If $A$ is a non-unital filtered $C^*$-algebra, we will by default endowed its unitalization $\tilde A$ with the filtration $\tilde A_E = A_E + \C$. A $*$-homomorphism $\phi : A \rightarrow B$ is said to be filtered if $\phi(A_E)\subseteq B_E$ for all $E\in\mathcal E$.\\

The crucial example for us will be crossed products of $G$-algebras by an étale groupoid $G$. Note that this definition generalizes that of \cite{OY2}. Indeed, as will be recalled later, the Roe algebras can be expressed as a crossed product by the so-called coarse groupoid, which is étale \cite{SkTuYu}, and the definition given here, applied to this groupoid, gives back a filtration equivalent to that of \cite{OY2}. Our filtration however does not depend on any choice of metric, but rather on the coarse class of the space. The second example will be that of crossed products by discrete quantum groups in the sense of Woronowicz.\cite{Wo}

%%%%%%%%%%%%%%%%%%%%%%%%%%%%%%%%%%%%%%%%%%%%%%%%%%%%%
\subsection{Almost unitaries and almost projections}
%%%%%%%%%%%%%%%%%%%%%%%%%%%%%%%%%%%%%%%%%%%%%%%%%%%%%

\begin{definition}
Let $(A,\mathcal E)$ be a unital filtered $C^*$-algebra. Let $\varepsilon\in(0,\frac{1}{4})$ and $E\in \mathcal E$ a controlled subset. The set of $\varepsilon$-$E$-unitaries is the set 
\[U^{\varepsilon, E}(A)= \{u\in A_E \text{ s.t. } ||u^*u-1||<\varepsilon\text{ and }||uu^*-1||<\varepsilon \}\]
and the set $\varepsilon$-$E$-projections is the set 
\[P^{\varepsilon, E}(A)= \{p\in A_E \text{ s.t. } p=p^*\text{ and }||p^2-p||<\varepsilon \}.\]
We will use the notation $P_n^{\varepsilon, E}(A)$ for $P^{\varepsilon, E}(M_n(A))$, and $U_n^{\varepsilon, E}(A)$ for $U^{\varepsilon, E}(M_n(A))$. Also, $P_\infty^{\varepsilon, E}(A)$ is the algebraic inductive limit of the $P_n^{\varepsilon, E}(A)$ under the natural inclusions
\[\left\{\begin{array}{rcl}
	P^{\varepsilon,E}_n(A) 		& \rightarrow	& P^{\varepsilon,E}_{n+1}(A)\\ 
	p 		& \mapsto 	& \begin{pmatrix}p& 0 \\ 0&0 \end{pmatrix}
\end{array}\right.\]
and $U_\infty^{\varepsilon, E}(A)$ is the algebraic inductive limit of the $U_n^{\varepsilon, E}(A)$ under the natural inclusions
\[\left\{\begin{array}{rcl}
	U^{\varepsilon,E}_n(A) 		& \rightarrow	& U^{\varepsilon,E}_{n+1}(A)\\ 
	u 		& \mapsto 	& \begin{pmatrix}u & 0 \\ 0& 1 \end{pmatrix}
\end{array}\right. .\]
\end{definition}

\begin{rk}Let $\varepsilon\in (0,\frac{1}{4})$ and $E\in\mathcal E$.\\
\begin{itemize}
\item[$\bullet$] If $p\in P^{\varepsilon,E}(A)$, then $p$ has a spectral gap around $\frac{1}{2}$, and functional calculus allows to define a genuine projection $\kappa_0(p)$, as in \cite{OY2}, by taking $\kappa_0$ to be a continuous function that vanishes inside the spectral gap and that is respectively $0$ and $1$ on the left and right part of the spectrum of $p$. For example, $\kappa_0=\chi_{(\frac{1}{2};\infty)}$ works. 
\item[$\bullet$] If $u\in U^{\varepsilon,E}(A)$, then $u^* u$ is invertible, and $u(u^* u)^{-\frac{1}{2}}$ defines a unitary, that we will denote $\kappa_1(u)$.
\end{itemize}
\end{rk}

In order to define controlled $K$-groups, define the following equivalence relations on $P^{\varepsilon, E}_\infty(A)\times \N$ and $U^{\varepsilon,E}_n(A)$.\\

\begin{itemize}

\item[$\bullet$] $(p,l) \sim (q,l')$ if there exists a homotopy of almost projections $h\in P^{\varepsilon, E}_\infty(A[0,1])$ and an integer $k$ such that 
\[h(0)=\begin{pmatrix} p & 0 \\ 0 & 1_{k+l'} \end{pmatrix} \text{ and }
h(1)=\begin{pmatrix} q & 0 \\ 0 & 1_{k+l} \end{pmatrix}\]
\item[$\bullet$] $u \sim v$ if there exists a homotopy of almost unitaries $h\in U^{3\varepsilon, E\circ E}_\infty(A[0,1])$ and an integer $k$ such that $h(0)= u \text{ and }h(1)=v$.\\
\end{itemize}

The following fact will be useful for future purposes. The reader can look at Proposition $1.30$ of \cite{OY3} for references. Recall the following definition.

\begin{definition}
Let $C>0$ and $A$ be a $C^*$-algebra. A map $h : [0,1]\rightarrow A $ is called $C$-Lipschitz if $||h(s)-h(t)||\leq C|s-t|$ for all $s,t\in [0,1]$.
\end{definition}

\begin{prop}\label{Lip}
There exists a universal constant $L>0$ such that, for any unital filtered $C^*$-algebra $(A,\mathcal E)$, any $\varepsilon\in(0,\frac{1}{4})$ and  any $E\in \mathcal E$, if $u_0$ and $u_1$ are homotopic in $U_n^{\varepsilon, E}(A)$, then there exists an integer $k$ and a $L$-lipschitz homotopy in $U_{n+k}^{3\varepsilon,E\circ E}$ connecting $\begin{pmatrix} u_0 & 0 \\ 0 & 1_k\end{pmatrix}$ and $\begin{pmatrix} u_1 & 0 \\ 0 & 1_k\end{pmatrix}$. 
\end{prop}

Denote $[(p,l)]_{\varepsilon,E}$ and $[u]_{\varepsilon,E}$ for the equivalence classes of almost-projections and almost-unitaries. Then, the same proof as \cite{OY2} shows that $[p,l]_{\varepsilon,E}+[q,l']_{\varepsilon,E}=[diag(p,q),l+l']_{\varepsilon,E}$ and $[u]_{\varepsilon,E}+[v]_{\varepsilon,E}=[diag(u,v)]_{\varepsilon,E}$ induces a group structure on the equivalence classes, that we denote $K_0^{\varepsilon,E}(A) = P^{\varepsilon, E}_\infty(A)\times \N / \sim$ and $K_1^{\varepsilon,E}(A) = U^{\varepsilon, E}_\infty(A) / \sim$.\\

If $A$ is not unital, let $\tilde A$ be the smallest unitalization of $A$, and $\rho_A: \tilde A \rightarrow \C; (a,\lambda)\mapsto \lambda$ the augmentation map. Then $K_0^{\varepsilon,E}(A)$ is defined as
\[\{[p,l]_{\varepsilon,E} : p\in P^{\varepsilon,E}_\infty (\tilde A), l\in \N \text{ s.t. rank}(\kappa_0(\rho_A(p)))=l \}\]
and $K_1^{\varepsilon,E}(A)$ is defined as $U_\infty^{\varepsilon,E}(\tilde A)/ \sim_{\varepsilon,E}$.\\

\begin{definition}
The controlled $K$-theory of a filtered $C^*$-algebra $(A,\mathcal E)$ is the family of abelian groups $\hat K_0(A) = (K_0^{\varepsilon,E}(A))_{\varepsilon\in (0,\frac{1}{4}),E\in\mathcal E}$ and $\hat K_1(A) = (K_1^{\varepsilon,E}(A))_{\varepsilon\in (0,\frac{1}{4}),E\in\mathcal E}$ defined above.\\
\end{definition}

We define canonical morphisms : if $\varepsilon, \varepsilon'\in (0,\frac{1}{4})$ and $E,E'\in\mathcal E$ such that $\varepsilon \leq \varepsilon'$ and $E \subseteq E'$, the natural homomorphism $K_*^{\varepsilon,E}(A)\hookrightarrow K_*^{\varepsilon',E'}(A)$ is denoted by $\iota_{\varepsilon,E}^{\varepsilon',E'}$. Notice that $\iota_{\varepsilon',E'}^{\varepsilon'',E''}\circ\iota_{\varepsilon,E}^{\varepsilon',E'}=\iota_{\varepsilon,E}^{\varepsilon'',E''}$ when this expression makes sense.\\

One has also forgetful morphisms $\iota_{\varepsilon,E} : K_*^{\varepsilon,E}(A)\rightarrow K_*(A)$ given by $[p,l]_{\varepsilon,E} \mapsto [\kappa_0(p)]-[1_l]$ and $[u]_{\varepsilon,E} \mapsto [\kappa_1(u)] $, and 
$\iota_{\varepsilon',E'}\circ\iota_{\varepsilon,E}^{\varepsilon',E'}=\iota_{\varepsilon, E}$ holds. The controlled $K$-theory groups approximate the usual $K$-groups in the sense of the following remarks. The reader is referred to \cite{OY2}, remark $1.17$, for a proof.

\begin{rk}
For every filtered $C^*$-algebra $(A,\mathcal E)$, any $y\in K(A)$ and any $\varepsilon\in (0,\frac{1}{4})$, there exists $E\in\mathcal E$ and $x\in K^{\varepsilon,E}(A)$ such that $\iota_{\varepsilon,E}(x) = y$.
\end{rk}

\begin{rk}\label{approximation}
There exists a universal constant $\lambda \geq 1$ such that, for every filtered $C^*$-algebra $(A,\mathcal E)$, every $\varepsilon\in (0,\frac{1}{4\lambda})$ and every $E\in \mathcal E$, the following holds :

Let $x\in K^{\varepsilon,E}(A)$ such that $ \iota_{\varepsilon,E}(x) = 0 $ in $K(A)$. Then, there exists $E'\in\mathcal E$ such that $E\leq E'$ and $\iota_{\varepsilon,E}^{\lambda\varepsilon,E'}(x)=0$ in $K^{\lambda\varepsilon,E'}(A)$.
\end{rk}

We list some examples that we will use, and some others that shall hopefully be developed in future work.

\begin{Expl}
Let $(X,\mathcal E)$ be a coarse space. The set of symmetric controlled subsets $\mathcal E_X$ is our proeminent example. Recall that it is the set of subsets $E\subseteq X\times X$ such that $\sup \{ d(x,y) : (x,y)\in E\}$ is finite and $E=E^{-1}$. Let $E,E'$ be such sets, their composition is given by :
\[E\circ E' = EE' \cup E'E \text{ where }EE' = \{(x,y) \text{ s.t. }\exists z\in X / (x,z)\in E \text{ and }(z,y)\in E'\}.\]
\end{Expl}

\begin{Expl}
Let $G$ be a $\sigma$-compact étale groupoid. Then the set of symmetric compact subsets $\mathcal E_G$ of $G$ is a coarse structure. It is the set of compact subsets $E\subseteq G$ such that $E=E^{-1}$, where $E^{-1} = \{g^{-1} : g\in E\}$. For $E$ and $E'$ in $\mathcal E_G$, their composition is given by :
\[E\circ E' = EE' \cup E'E \text{ where }EE' = \{gg' : (g,g')\in G^{(2)}, g\in E\text{ and } g'\in E'\}.\] 
If $G$ is $\sigma$-compact, and $A$ is a $G$-algebra, $A\rtimes_r G$ is naturally filtered by $\mathcal E_G$ : if $E\subseteq G$ is a compact subset, define $(A\rtimes_r G)_E = \{f\in C_c(G,A) : \text{ supp}(f)\subseteq E\}$. This situation has important particular cases :
\begin{itemize}
\item[$\bullet$] Let $G$ be the coarse groupoid of a coarse space $(X,\mathcal E)$, which is étale. Then $\mathcal E_G$ is given by \[\{\overline{E} : E\in\mathcal E_X\} \cong \mathcal E_X.\] 
\item[$\bullet$] Let $\Gamma$ be a finitely generated group acting by homeomorphism on a topological space $X$. Recall that the word length $l$ defines a proper metric on $\Gamma$. Define, for $R>0$ and $K\subseteq X$ compact, 
\[\Delta_{R,K} = \{(x,g)\in G \text{ s.t. } l(g)\leq R \text{ and } x\in K\}\] 
and $\mathcal E_G $ as the set of $E\subseteq G$ such that $\exists R>0, E\subseteq \Delta_R$ and $E = E^{-1}$. It provides the étale action groupoid $G=X\rtimes \Gamma$ with a coarse structure.
\end{itemize}
\end{Expl}

\begin{Expl}\label{EQG}
Let $\mathbb G$ be a compact quantum group in the sense of Woronowicz. Denote by $\hat{\mathbb G}$ its discrete dual. A unitary representation is said to be symmetric if it is equivalent to its contragredient. Then, the set $\mathcal E_{\mathbb G}$ of symmetric unitary representations of $\mathbb G$ is a coarse structure w.r.t. the composition $\pi\circ \pi'=(\pi\otimes \pi')\oplus (\pi\otimes \pi')$, and $\pi\leq\pi'$ if $\pi$ is equivalent to a subrepresentation of $\pi'$. Moreover, for every $\hat{\mathbb G}$-algebra $A$,  $A\rtimes_r \hat{\mathbb G}$ is $\mathcal E_{\mathbb G}$-filtered.
\end{Expl}

%%%%%%%%%%%%%%%%%%%%%%%%%%%%%%%%%%
\subsection{Quantitative objects}
%%%%%%%%%%%%%%%%%%%%%%%%%%%%%%%%%%

In order to study functorial properties of controlled $K$-theory, we will adapt and study the notion of quantitative object defined in \cite{OY2}.

\begin{definition}
A quantitative object is a family of abelian groups $\hat{\mathcal{O}}=\{\mathcal{O}_{\varepsilon, E}\}_{\varepsilon\in (0,\frac{1}{4}), E\in\mathcal{E}}$ endowed with a family of group homomorphisms $\phi_{\varepsilon, E}^{\varepsilon', E'} : \mathcal{O}_E\rightarrow \mathcal{O}_{E'}$ for any $E,E'\in\mathcal E$ and $0<\varepsilon\leq \varepsilon'<\frac{1}{4}$ such that $E\subseteq E'$, satisfying $\phi_{\varepsilon, E}^{\varepsilon, E}= id_{\mathcal{O}_E}$ and $\phi_{\varepsilon', E'}^{\varepsilon'', E''}\circ \phi_{\varepsilon, E}^{\varepsilon', E'} =\phi^{\varepsilon'', E''}_{\varepsilon, E}$ whenever $E\leq E' \leq E''$ and $\varepsilon<\varepsilon'<\varepsilon''$.
\end{definition}

We need to define controlled morphisms between quantitative objects. We first define control pairs, which are essentially what ensures that the controlled morphisms do not distort too much the propagation.\\

\begin{definition}
A control pair is a couple $\rho=(a,h)$ where $a\in (0,\frac{1}{4})$ and $h : (0,\frac{1}{4a})\rightarrow \N^*$ is a non-decreasing function. 
\end{definition}

Control pairs can be naturally composed, and if $(a,h)$ and $(b,k)$ are two control pairs, then their composition, denoted by $(b,k)\ast(a,h)$, is defined by $(ab,k\ast h)$, where $k\ast h : (0,\frac{1}{4ab})\rightarrow \N^* ; \varepsilon \mapsto k_{a \varepsilon}h_\varepsilon$. \\

Control pairs naturally act on the index subset of the controlled $K$-groups. Indeed, if $\varepsilon\in (0,\frac{1}{4a})$ and $E\in\mathcal E$, $(a,h).(\varepsilon,E)= (a\varepsilon,E^{h_\varepsilon})$ is in $(0,\frac{1}{4})\times\mathcal E$. This allows to define controlled morphims.\\

We can also compare control pairs. Indeed, define the following partial order : $(a,h) \leq (b,k)$ if $a \leq b$ and $h_\varepsilon\leq k_\varepsilon$ for all $\varepsilon\in (0,\frac{1}{4a})$.\\ 

\begin{definition}
Let $\hat{\mathcal O}$ and $\hat{\mathcal O'}$ be two quantitative objects and $\rho=(a,h)$ a control pair. A $\rho$-controlled morphism is a family of groups homomorphims $F_{\varepsilon, E} : \mathcal O_{\varepsilon, E} \rightarrow \mathcal{ O'}_{a\varepsilon, E^{h_\varepsilon}}$ for any $\varepsilon\in(0,\frac{1}{4a})$ and $E\in\mathcal E$, such that
\[\phi^{a\varepsilon',E'^{h_\varepsilon'}}_{a\varepsilon,E^{h_\varepsilon}} \circ F_{\varepsilon, E} =F_{\varepsilon',E'} \circ \phi_{\varepsilon,E}^{\varepsilon',E'}\] for any $0<\varepsilon<\varepsilon'<\frac{1}{4a}$ and $E\subseteq E'$.
\end{definition}

\begin{rk}
%If we don't specify any control pair, it is implicit and evident from the context, or it is not crucial to specify it. 
When not specified, the control pair is evident from the context. For example, we will often refer to a controlled morphism, meaning a $\alpha$-controlled morphism for some control pair $\alpha$. For a controlled morphism $\hat F : \hat K(A)\rightarrow \hat K(B)$, we will denote $F:K(A)\rightarrow K(B)$ the unique homomorphism it induces in $K$-theory. We will always try to indicate an analogy with the classical case (as opposed to the controlled or quantitative case) by putting a hat on top of controlled objects that are hopefully inducing a well known object. For example, controlled $K$-theory is $\hat K$, the controlled assembly map will be denoted $\hat \mu_{G,A}$, etc.\\
\end{rk}

Let $\rho=(\lambda,h),\alpha,\beta$ be control pairs, and $F :  \hat{\mathcal O} \rightarrow \hat{\mathcal O'}$ and $G : \hat{\mathcal O} \rightarrow \hat{\mathcal O'}$ be $\alpha$- and $\beta$-controlled morphisms respectively. We write $F\sim_\rho G$ if :
\begin{itemize}%{$bullets$}
\item[$\bullet$] $\alpha \leq \rho$ and $\beta \leq \rho$,
\item[$\bullet$] for any $\varepsilon\in (0,\frac{1}{4\lambda})$ and $E\in \mathcal E$, the following diagram commutes : 
\[\begin{tikzcd}
 \mathcal O_{\varepsilon,E} \arrow{r}{F_{\varepsilon,E}} \arrow{d}{G_{\varepsilon,E}} & \mathcal O'_{\alpha (\varepsilon,E) } 
\arrow{d}{\iota_{\beta(\varepsilon,E)}^{\rho(\varepsilon,E)}} \\
 \mathcal O'_{\beta (\varepsilon,E) } \arrow{r}{\iota_{\alpha(\varepsilon,E)}^{\rho(\varepsilon,E)}} & \mathcal O'_{\rho(\varepsilon,E)}\\
\end{tikzcd}.\]
\end{itemize}

\begin{rk} \label{rk2.5} The following statement can be found in \cite{OY2} (remark $2.5$). Let $F,F' : \mathcal O\rightarrow \mathcal O'$, $G: \mathcal O'\rightarrow \mathcal O''$ and $G': \mathcal O'''\rightarrow \mathcal O$ be $\alpha$-,$\alpha'$-, $\beta$- and $\beta'$-controlled morphisms respectively. Let $\rho $ be a control pair such that $F\sim_\rho F'$. Then $G\circ F \sim_{\beta\ast\rho}  G\circ F'$ and $F\circ G'\sim_{\rho \ast \beta'} F'\circ G'$.
\end{rk}

\begin{definition}
Let $\alpha$ and $\rho$ be control pairs satisfying $\alpha \leq \rho$ and $F : \hat{\mathcal O} \rightarrow \hat{\mathcal O'}$ a $\alpha$-controlled morphism. We say that $F$ is $\rho$-invertible if there exist a controlled morphism $ G : \hat{\mathcal O'} \rightarrow \hat{\mathcal O}$ such that $G \circ F \sim_\rho Id_{\hat{\mathcal O}}$ and $F\circ G \sim_\rho Id_{\hat{\mathcal O'}}$. $G$ is said to be a $\rho$-inverse for $F$.
\end{definition}

As we will see for controlled assembly maps, the correct notions of injectivity and surjectivity for controlled morphisms need to be adpated in the following way.\\

\begin{definition}
Let $\rho=(\lambda,h)$ and $\alpha$ be controlled pairs, and $F : \hat{\mathcal O} \rightarrow \hat{\mathcal O'}$ a $\alpha$-controlled morphism. 
\begin{itemize}
\item[$\bullet$] $F$ is $\rho$-injective if, given any $\varepsilon \in (0,\frac{1}{4\lambda})$ and $E\in \mathcal E$, $\alpha\leq \rho$ and, for all $x\in \mathcal O_{\varepsilon, E}$ such that $F_{\varepsilon, E}(x)=0$, then $\iota_{\varepsilon, E}^{\lambda\varepsilon, h_\varepsilon E}(x)=0$,
\item[$\bullet$] $F$ is $\rho$-surjective if, given any $\varepsilon \in (0,\frac{1}{4\lambda})$ and $E\in \mathcal E$, for any $y\in \mathcal O' _{\varepsilon, E}$, there exists $x\in \mathcal O_{\rho(\varepsilon,E)}$ such that $F_{\rho(\varepsilon,E)}(x)= \iota_{\varepsilon,E}^{\rho(\varepsilon,E)}$,
\end{itemize}
\end{definition}

\begin{rk}
If $F$ is a $\rho$-isomorphism, then there exists a control pair $\rho'$ only depending on $\rho$ such that $F$ is both $\rho'$-injective and $\rho'$-surjective.
\end{rk}
%%%%%%%%%%%%%%%%%%%%%%%%%%%%%%%%%%%%%%%
\subsection{Controlled exact sequences}
%%%%%%%%%%%%%%%%%%%%%%%%%%%%%%%%%%%%%%%

\begin{definition}
Let $F : \hat{\mathcal O}\rightarrow \hat{\mathcal O'}$ and $G : \hat{\mathcal O'}\rightarrow \hat{\mathcal O''}$ be $(\alpha,h)$-controlled and a $(\beta,k)$-controlled morphisms respectively. The sequence
\[\begin{tikzcd}[column sep = small] \hat{\mathcal O} \arrow{r}{F} & \hat{\mathcal O'} \arrow{r}{G} & \hat{\mathcal O''} \end{tikzcd}\]
is called $\rho$-exact at $\hat{\mathcal O'}$ if $G\circ F=0$ and if for all $\varepsilon\in (0,\frac{1}{4 \max (\lambda \alpha,\beta)})$, $E\in\mathcal E$ and any $y\in \mathcal O'_{\varepsilon,E}$ such that $G_{\varepsilon,E}(y) = 0$, then there exists $x\in \mathcal O_{\rho(\varepsilon,E)}$ such that $F_{\rho(\varepsilon,E)}(x)=\iota_{\varepsilon,E}^{\rho(\varepsilon,E)} (y)$.\\
A sequence of controlled morphisms 
\[\begin{tikzcd}[column sep = small] ... \arrow{r}{F_{k-2}} & \hat{\mathcal O_{k-2}} \arrow{r}{F_{k-1}} & \hat{\mathcal O_{k-1}} \arrow{r}{F_k} & \hat{\mathcal O_k} \arrow{r}{F_{k+1}} & ... \end{tikzcd}\] 
is said to be $\rho$-exact if the sequence
\[\begin{tikzcd}[column sep = small] \hat{\mathcal O_{j-1}} \arrow{r}{F_{j}} & \hat{\mathcal O_{j}} \arrow{r}{F_{j+1}} & \hat{\mathcal O_{j+1}} \end{tikzcd}\]
is $\rho$-exact at $\hat{\mathcal O_j}$ for all $j$.
\end{definition}

%%%%%%%%%%%%%%%%%%%%%%%%%%%%%%%%%
\subsection{Morita equivalence}
%%%%%%%%%%%%%%%%%%%%%%%%%%%%%%%%%

A controlled version of the Morita equivalence exists. Indeed, the classical Morita equivalence states that, if $e$ is a rank $1$ projection in $\mathcal L (H)$, the map $A\rightarrow A\otimes \mathfrak K(H) ; a\mapsto a\otimes e$ induces an isomorphism in $K$-theory. But this map preserves propagation, hence the same proof as in \cite{OY2} gives the following proposition.

\begin{prop}
Let $A$ be a $\mathcal E$-filtered $C^*$-algebra and $H$ a separable Hilbert space. Then the $*$-homomorphism
\[A\rightarrow A\otimes \mathfrak K(H) ; \quad a\mapsto 
\begin{pmatrix}a & & \\  & 0 & \\ & & ... \end{pmatrix}\]
induces a group isomorphism 
\[\mathcal M_A^{\varepsilon,E} : K^{\varepsilon,E}(A)\rightarrow K^{\varepsilon,E}(A\otimes \mathfrak K(H)) \]
for every $\varepsilon\in(0,\frac{1}{4})$ and $E\in\mathcal E$. The family $\mathcal M_A = (\mathcal M_A^{\varepsilon,E} )_{\varepsilon\in(0,\frac{1}{4}),E\in\mathcal E}$ is called the controlled Morita equivalence and is a controlled morphism. It induces the usual Morita equivalence $M_A: K(A)\rightarrow K(A\otimes \mathfrak K(H))$ in $K$-theory. 
\end{prop}

%%%%%%%%%%%%%%%%%%%%%%%%%%%%%%%%%%%%%%%%%%%%%%%%%
\subsection{Controlled $6$-term exact sequences}
%%%%%%%%%%%%%%%%%%%%%%%%%%%%%%%%%%%%%%%%%%%%%%%%%

We will describe the $6$-term controlled exact sequence associated to a completely filtered extension of $C^*$-algebras. For any extension of $C^*$-algebras 
\[\begin{tikzcd}[column sep = small]0\arrow{r} & J\arrow{r} & A\arrow{r} & A/J\arrow{r} & 0 \end{tikzcd},\]
we denote $\partial_{J,A}$ the boundary map $K_*(A/J)\rightarrow K_{*+1}(J)$. The reader can find all the proofs and properties of the following results in \cite{OY2}.\\

We fix a coarse structure $\mathcal E$, and we will consider $\mathcal E$-filtered $C^*$-algebras.

\begin{definition}
Let $A$ a filtered $C^*$-algebra and $J\subseteq A$ an ideal. If $J_E = A_E\cap J$, the extension
\[\begin{tikzcd}[column sep = small]0\arrow{r} & J\arrow{r} & A\arrow{r} & A/J\arrow{r} & 0 \end{tikzcd}\]
is said to be completely filtered if the continuous linear bijection $A_E/J_E \hookrightarrow (A_E+J)/J$ induced by the inclusion $A_E\hookrightarrow A$ is a complete isometry, i.e.
\[ \inf_{y\in M_n(J_E)} ||x+y|| = \inf_{y\in M_n(J)} ||x+y||\quad,\forall n\in \N,x\in M_n(A_E),E\in \mathcal E.\]
\end{definition}

\begin{prop}
There exists a control pair $(\alpha_D,k_D)$ such that for any completely filtered extension of $C^*$-algebras
\[\begin{tikzcd}[column sep = small]0\arrow{r} & J\arrow{r} & A\arrow{r} & A/J\arrow{r} & 0 \end{tikzcd}\]
there exists a $(\alpha_D,k_D)$-controlled morphism of odd degree
\[D_{J,A} : \hat K(A/J)\rightarrow \hat K(J)\]
which induces $\partial_{J,A}$ in $K$-theory.
\end{prop}

We will denote 
\begin{itemize}
\item[$\bullet$] $D_A$ the controlled boundary map associated to the completely filtered extension of $C^*$-algebras $0 \rightarrow SA \rightarrow CA \rightarrow A \rightarrow 0 $,
\item[$\bullet$] $D^j_{J,A}$ the restriction of $D_{J,A}$ to $\hat K_j(A/J)\rightarrow \hat K_{j+1}(J)$.
\end{itemize} 

\begin{thm}
There exists a universal control pair $(\lambda,h)$, which does not depend on $\mathcal E$, such that for any completely filtered extension of $C^*$-algebras 
\[\begin{tikzcd}[column sep = small]0\arrow{r} & J\arrow{r}{\iota} & A\arrow{r}{q} & A/J\arrow{r} & 0 \end{tikzcd}\]
the following $6$-term exact sequence is $(\lambda,h)$-exact
\[\begin{tikzcd}[column sep = small]
 \hat K_1(J) \arrow{r}{\iota_*} & \hat K_1(A) \arrow{r}{q_*} & \hat K_1(A/J)\arrow{d}{D^1_{J,A}} \\
 \hat K_0(A/J) \arrow{u}{D^0_{J,A}} & \hat K_0(A) \arrow{l}{q_*} & \hat K_0(J)\arrow{l}{\iota_*}
\end{tikzcd}.\]
\end{thm}

The following remark can be found in \cite{OY2} (remark $3.8$) and will be used to prove functorial properties of the controlled Roe and Kasparov transformations.

\begin{rk}\label{rk3.8}
Let $A$ and $B$ two $\mathcal E$-filtered $C^*$-algebras, and $\phi : A\rightarrow B$ a filtered $*$-homomorphism. Let $I$ and $J$ be respectively ideals in $A$ and $B$ and assume that :
\begin{itemize}
\item[$\bullet$] $0 \rightarrow I \rightarrow A \rightarrow A/ I \rightarrow 0$ and $0 \rightarrow J \rightarrow B \rightarrow B/J \rightarrow 0$ are completely filtered extensions of $C^*$-algebras,
\item[$\bullet$] $\phi(I)\subseteq J$,
\end{itemize}
then $D_{J,B}\circ \tilde\phi_* = \phi_* \circ D_{I,A}$.
\end{rk}

%%%%%%%%%%%%%%%%%%%%%%%%%%%%%%%%%%%%%%%%%%
\subsection{Tensorisation in $KK$-theory}
%%%%%%%%%%%%%%%%%%%%%%%%%%%%%%%%%%%%%%%%%%

If $B$ is a filtered $C^*$-algebra and $A$ any $C^*$-algebra, and if $A\otimes B$ is the spatial tensor product, define, for all $E\in \mathcal E$, $(A\otimes B)_E$ as the closure of the algebraic tensor product $A\otimes B_E$ in $A\otimes B$. Then $(A\otimes B_E)_{E\in\mathcal E}$ defines a filtration of $A\otimes B$. If $\phi : A_1 \rightarrow A_2$ is a $*$-homomorphism, we use the notation $\phi_B$ for the induced $*$-homomorphism $A_1\otimes B\rightarrow A_2\otimes B$. \\

In \cite{kasparovKKNovikov}, G. Kasparov defined a map
\[\tau_B : KK(A_1,A_2)\rightarrow KK(A_1\otimes B, A_2\otimes B)\]
for any $C^*$-algebras $A_1$ and $A_2$, which is compatible with the Kasparov product. Any $z\in KK(A_1,A_2)$ defines a morphism
\[K(A_1\otimes B)\rightarrow K(A_2\otimes B)\]
which is proved in \cite{OY2} to be induced from a controlled morphism. The following theorem is borrowed from \cite{OY2}.\\

\begin{thm}\label{tensorization}
There exists a control pair $(\alpha_\tau,k_\tau)$ such that, for any filtered $C^*$-algebra $B$, any $C^*$-algebras $A_1$ and $A_2$ and any $K$-cycle $z\in KK(A_1,A_2)$, there exists a $(\alpha_\tau,k_\tau)$-controlled morphism $\hat \tau_B : \hat K(A_1\otimes B)\rightarrow \hat K(A_2\otimes B)$
such that :
\begin{itemize}
\item[$\bullet$] $\hat \tau_B(z)$ induces right-multiplication by $\tau_B(z)$ in $K$-theory,
\item[$\bullet$] for any $K$-cycles $z,z'\in KK(A_1,A_2)$, $\hat \tau_B(z+z')=\hat\tau_B(z)+\hat\tau_B(z')$,
\item[$\bullet$] if $\phi : A_1\rightarrow A'_1$ is a $*$-homomorphism, then $\hat\tau_B(\phi^*(z)) =  \hat\tau_B(z)\circ (\phi_B)_*$ for any $z\in KK(A'_1,A_2)$,
\item[$\bullet$] if $\phi : A_2'\rightarrow A_2$ is a $*$-homomorphism, then $\hat\tau_B(\phi_*(z)) = (\phi_B)_*\circ \hat\tau_B(z)$ for any $z\in KK(A_1,A'_2)$,
\item[$\bullet$] $\hat \tau_B([Id_A])\sim_{(\alpha_\tau,k_\tau)} Id_{\hat K(A\otimes B)}$,
\item[$\bullet$] for any $C^*$-algebra $D$, any $K$-cycle $z\in KK(A_1,A_2)$, $\hat\tau_B (\tau_D(z))= \hat\tau_{B\otimes D}(z)$.
\item[$\bullet$] for any semi-split extension $\begin{tikzcd}[column sep = small] 0 \arrow{r}& J \arrow{r} & A \arrow{r} & A/J\arrow{r} & 0\end{tikzcd}$ with boundary element $[\partial_{J,A}]\in KK_1(A/J,J)$, $\hat\tau_B([\partial_{J,A}])=D_{J\otimes B,A\otimes B}$.
\end{itemize}
\end{thm}

This controlled tensorisation map respects Kasparov product. See \cite{OY2} for a proof.

\begin{thm}
There exists a control pair $\lambda$ such that, for any separable $C^*$-algebras $A_1$ and $A_2$, any filtered $C^*$-algebra $B$, the following holds : for any $z\in KK(A_1,A_2)$ and $z'\in KK(A_2,A_3)$,
\[\hat\tau_B(z\otimes z')\sim_\lambda \hat\tau_B(z')\circ\hat\tau_B(z)\]
\end{thm}

%%%%%%%%%%%%%%%%%%%%%%%%%%%%%%%%%%%%%%%%%%
\subsection{Controlled Bott periodicity}
%%%%%%%%%%%%%%%%%%%%%%%%%%%%%%%%%%%%%%%%%%

We recall in this section how to construct controlled Bott maps. \\

Recall the definition of the Toeplitz extension. Let $v\in \mathcal L(l^2(\N))$ be the unilateral shift on $l^2(\N)$, i.e. $v(e_n) = e_{n+1}$ if $(e_n)_{n\in \N}$ is the canonical basis. Then $v$ is an isometry, i.e. $v^*v = 1$. Let us denote $\mathcal T $ the $C^*$-algebra of $\mathcal L(l^2(\N))$ generated by $v$. Then the sequence of $C^*$-algebras 
\[0\rightarrow \mathfrak K(l^2(\N)) \rightarrow \mathcal T \rightarrow C(\mathbb S) \rightarrow 0 \]
is exact. Set $\mathcal T_0$ the preimage of $C_0(0,1)$ under the last arrow, so that the sequence of $C^*$-algebras
\[0\rightarrow \mathfrak K(l^2(\N)) \rightarrow \mathcal T_0 \rightarrow C_0(0,1) \rightarrow 0 \]
is exact. This last sequence is called the Toeplitz extension, let $T\in KK_1(S, \mathfrak K(l^2(\N)) )$ denotes its boundary.\\

Let us put $S = C_0(0,1)$ and $C= C_0(0,1]$ , so that evaluation at $1$ induces the following exact sequence of $C^*$-algebras $0\rightarrow S \rightarrow C \rightarrow \C \rightarrow  0 $. Let us denote by $[\partial]\in KK_1(\C,S)$ the class in $KK$-theory of this extension. It is called the Bott element, and is invertible in $KK_1(\C,S)$, its inverse being the boundary of the Toeplitz extension, up to Morita equivalence (see \cite{OY2}, section $4.2$). Recall that the boundary of this extension is given by right Kasparov product by $[\partial]$.\\ 

For any $C^*$-algebra $A$, we can tensorize the previous extension by $A$ to have an exact sequence $0\rightarrow SA \rightarrow CA \rightarrow A \rightarrow  0 $. Its class in $KK$-theory, denoted $[\partial_{A}]\in KK_1(A,SA)$, actually satisfies $\tau_A([\partial ]) =[\partial_{A}]$. Let us denote $\partial_A : K_*(A)\rightarrow K_*(SA)$ the odd degree boundary homomorphism of the extension, it is induced by right multiplication by $[\partial_A]$. It is also invertible, as we can tensorize the Toeplitz extension to get the following exact sequence of $C^*$-algebras
\[0\rightarrow \mathfrak K(l^2(\N))\otimes A \rightarrow \mathcal T_0 \otimes A \rightarrow SA \rightarrow 0 \]
whose boundary element $T_A =[\partial_{\mathcal T_0 \otimes A }]\in KK_1(SA,\mathfrak K(l^2(\N))\otimes A)$ satisfies $\tau_A(T)= T_A$. As $\tau_A$ respects Kasparov product, one can see that $[\partial_A]$ and $T_A$ are mutual inverse up to Morita equivalence.\\

Tensorization by $S$ and $C$ is functorial with respect to $*$-homomorphisms, and if $g : A\rightarrow B$ is a $*$-homomorphism, we denote $Sg : SA\rightarrow SB$ and $Cg : CA\rightarrow CB$ the induced $*$-homomorphisms. By naturality of boundary maps w.r.t. extensions, we get the following relation :
\[(Sg)_*[\partial_{SA}] = [\partial_{SA}]\otimes_{SA} [Sg] = [g]\otimes_{B} [\partial_{SB}] = g^*[\partial_{SB}] .\]

%Recall that there exist controlled morphisms $D_{A} :\hat K_0(A)\rightarrow \hat K_{1}(SA)$ and $Z_A :\hat K_1(SA)\rightarrow \hat K_{0}(A\otimes \mathfrak K)$ that induce right-multiplication by $[\partial_A]$ and $T_A$ respectively, i.e.
%\[\iota_{\varepsilon,E}\circ D_A(x) = \kappa_0(x) \otimes [\partial_A]\quad \iota_{\varepsilon,E}\circ Z_A(y) = \kappa_1(y) \otimes T_A \quad \forall x\in K_0^{\varepsilon, E}(A),y\in K_1^{\varepsilon, E}(SA). \]

Recall that, if $A$ is unital and $p$ is a projection in $A$, $\partial_A [p]$ is given by the homotopy class of the path of unitaries 
\[\left\{\begin{array}{rcl} [0,1] & \rightarrow & U(\tilde{SA}) \\ t & \mapsto & e^{2i\pi t}p + 1-p\end{array}\right.\]
One can perform a similar construction in term of almost-projections. Let $q\in P_n^{\varepsilon,E}(A)$ and $l$ an integer, define 
\[V_{q,l}(t) = \text{diag }(e^{-2i\pi l t}, 1, 1,...,1) \ (e^{2i\pi t}q +1-q)\in U^{5\varepsilon,E}_n(\tilde{SA}),\]
then 
\[\left\{\begin{array}{rcl} 
K_0^{\varepsilon,E}(A) & \rightarrow & K_1^{\varepsilon,E}(SA) \\ 
\ [q,l]_{\varepsilon,E} & \mapsto & [V_{q,l}]_{5\varepsilon,E}
\end{array}\right.\]
defines a $(5,\textbf{1})$-controlled morphism $Z_A : \hat K_0(A)\rightarrow \hat K_1(SA)$, where $\textbf{1}$ is the constant function $\textbf{1}_\varepsilon(E) = 1$, $\forall E\in\mathcal E$.

\begin{prop}[\cite{OY2} Prop. $3.9$ and $3.10$]
Let $A$ be a unital $C^*$-algebra. There exist universal control pairs $(\lambda,h)$ and $(\lambda',h')$ such that 
\[\begin{array}{c} D^0_A \sim_{(\lambda,h)} Z_A \quad \text{and}\\
D^1_{ \mathfrak K(l^2(\N)) \otimes A,\mathcal T_0\otimes A}\circ Z_A\sim_{(\lambda',h')} \mathcal M_A \end{array}\]
where $\mathcal M_A$ is the controlled Morita equivalence.
\end{prop}

%%%%%%%%%%%%%%%%%%%%%%%%%%%

%%%%%%%%%%%%%%%%%%%%%%%%%%%
\section{Controlled assembly maps for étale groupoids}

In this section, we will always use a coarse structure $\mathcal E$ of a locally compact $\sigma$-compact étale groupoid generated by a countable subset of its compact symmetric subsets such that for every compact subset $K\subseteq G$, there exists a $E\in\mathcal E$ such that $K\subseteq E$.

%%%%%%%%%%%%%%%%%%
\subsection{Controlled Kasparov transform} %%
%%%%%%%%%%%%%%%%%%

%Let $A$ and $B$ be two $G$-$C^*$-algebras, and $H$ a separable Hilbert space, $l^2(\Z)$ for instance, and $H_G= H\otimes L^2(G,\lambda)$. The standard Hilbert module over $B$ is denoted by $H_B=H_G\otimes B$, and $K_B$ is the algebra of compact operators for $H_B$, i.e. $K(H)\otimes L^2(G,\lambda)\otimes B$. \\

In this section, for every $G$-algebra $B$, we use the notation $K_{B\rtimes_r G}$ for $\mathfrak K_{B \rtimes_r G}(B \rtimes_r G)\cong \mathfrak K \otimes (B\rtimes_r G)$.\\

Let $A$ and $B$ be two $G$-$C^*$-algebras, with actions denoted by $\alpha : s^* A \rightarrow r^* A$ and $\beta : s^* B \rightarrow r^* B$. %is the action of $G$, notice that the canonical action on the standard Hilbert module $H_B$ is given by $V=id_H\otimes \beta : s^* H_B \cong H_{s^* B} \rightarrow r^* H_B \cong H_{r^* B}$. \\
Recall that every $K$-cycle $z\in KK^G(A,B)$ can be represented as a triple $(H_B, \pi, T)$ where :
\begin{itemize}
\item[$\bullet$] $H_B$ is equipped with an action $V\in\mathcal L_{s^* B}(s^* H_B , r^* H_B)$ of $G$.
\item[$\bullet$]$\pi : A\rightarrow \mathcal L_B(H_B)$ is a $G$-equivariant $*$-homomorphism.
\item[$\bullet$]$T\in \mathcal L_B(H_B)$ is a self-adjoint operator.
\item[$\bullet$] $T$ and $\pi$ satisfy the $K$-cycle conditions, i.e. $[T,\pi(a)]$, $\pi(a)(T^2-id_{H_B})$ are compact operators in $\mathfrak K_B(H_B)$ and $\pi(a)(r^*T-Vs^*T V^*)$ are compact operators in $\mathfrak K_{r^* B}(r^* H_B)\cong \mathfrak K_{r^* B}(H_{r^* B})$ for all $a\in A, g\in G$.\\
\end{itemize}

Set $T_G= T\otimes id_{B\rtimes_r G}\in \mathcal L_{B\rtimes_r G}(H_B\otimes (B\rtimes_r G))\simeq \mathcal L_{B\rtimes_r G}(H_{B\rtimes_r G})$, and $\pi_G: A\rtimes_r G \rightarrow \mathcal L_{B\rtimes_r G}(H_{B\rtimes_r G})$. Then, according to Le Gall \cite{LeGall}, $(H_{B\rtimes_r G}, \pi_G, T_G)$ represents the $K$-cycle $j_G(z)\in KK(A\rtimes_r G,B\rtimes_r G)$. Let us construct a controlled morphism associated to $z$,
\[J_G(z) : \hat K(A\rtimes_r G)\rightarrow \hat K(B\rtimes_r G), \]
which induces right multiplication by $j_G(z)$ in $K$-theory.

\subsubsection{Odd case}

Let us start with an odd element $z\in KK_1^G(A,B)$. Let $(H_B,\pi,T)$ be a $K$-cycle representing $z$. Set $P=\frac{1+T}{2}$ and $P_G=P\otimes id_{B\rtimes_r G}$. We define
\[E^{(\pi,T)}=\{(x,P_G\pi_G(x)P_G + y) : x\in A\rtimes_r G, y\in K_{B\rtimes_r G}\}\]
a $C^*$-algebra which is filtered by
\[E_U^{(\pi,T)}=\{(x,P_G\pi_G(x)P_G + y) : x\in (A\rtimes G)_U, y\in K\otimes (B\rtimes G)_U\}\]
for all $U\in\mathcal E$. This $C^*$-algebra fits into the filtered extension
\[\begin{tikzcd}[column sep = small]
0\arrow{r} & K_{B\rtimes_r G}\arrow{r} & E^{(\pi,T)} \arrow{r} & A\rtimes_r G \arrow{r}& 0
\end{tikzcd}\]
which is semi split by  $s :\left\{\begin{array}{lll}A\rtimes_r G & \rightarrow & E^{(\pi,T)} \\ x & \mapsto & (x, P_G \pi_G(x)P_G)\end{array}\right.$.\\

Let us show that the controlled boundary map of this extension does not depend on the representant chosen, but only on the class $z$.

\begin{lem}\label{ClassIndepedance} With the above notations, the controlled boundary map $D_{K_{B\rtimes_r G},E^{(\pi,T)}}$ only depends on the class $z$.
\end{lem}

\begin{dem}
Let $(H_B, \pi_j,T_j), j=0,1$ two $K$-cycles which are homotopic via $(H_{B[0,1]},\pi,T)$. We denote $e_t$ the evaluation at $t\in[0,1]$ for an element of $B[0,1]$, and set $y_t=e_t(y)$ for such a $y$. The $*$-morphism
\[\phi : \left\{\begin{array}{lll}E^{(\pi,T)} & \rightarrow & E^{(\pi_t,T_t)} \\ (x,y) & \mapsto & (x, y_t)\end{array}\right.\]
satisfies $\phi(K_{B[0,1] \rtimes_r G})\subseteq K_{B \rtimes_r G}$ and makes the following diagram commute
\[\begin{tikzcd}[column sep = small]
0\arrow{r} & K_{B[0,1] \rtimes_r G}\arrow{r}\arrow{d}{\phi_{|K_{B[0,1] \rtimes_r G}}} & E^{(\pi,T)} \arrow{r}\arrow{d}{\phi} & A\rtimes_r G \arrow{r}\arrow{d}{=}& 0 \\
0\arrow{r} & K_{B \rtimes_r G}\arrow{r} &  E^{(\pi_t,T_t)} \arrow{r} & A\rtimes_r G \arrow{r} & 0
\end{tikzcd}.\]

According to remark \ref{rk3.8}, the following holds
\[D_{K_{B\rtimes_r G},E^{(\pi_t,T_t)}} = \phi_* \circ D_{K_{B[0,1]\rtimes_r G},E^{(\pi,T)}}.\]
As $id \otimes e_t$ gives a homotopy between $id\otimes e_0$ and $id\otimes e_1$, and as if two $*$-morphisms are homotopic, then they are equal in controlled $K$-theory, 
\[D_{K_{B\rtimes_r G}, E^{(\pi_0,T_0)}}=D_{K_{B\rtimes_r G}, E^{(\pi_1,T_1)}}\]
holds, and the boundary of the extension $E^{(\pi,T)}$ depends only on $z$.\\
\qed
\end{dem}

\begin{definition}
The controlled Kasparov transform of an element $z\in KK_1^G(A,B)$ is defined as the composition
\[J_{red,G}(z)=\mathcal M_{B\rtimes_r G}^{-1}\circ D_{K_{B\rtimes_r G}, E^{(\pi,T)}}.\]
\end{definition}

As the boundary map is a $(\alpha_D,k_D)$-controlled morphism and the Morita equivalence preserves the filtration, $J_{red,G(z)}$ is  $(\alpha_D,k_D)$-controlled. 

\begin{prop}\label{Kasparov1}
Let $A$ and $B$ two $G$-$C^*$-algebras. There exists a control pair $(\alpha_J,k_J)$ such that for every $z\in KK^G_1(A,B)$, there exists a $(\alpha_J,k_J)$-controlled morphism
\[J_{red,G}(z) : \hat K_*(A\rtimes_r G)\rightarrow \hat K_{*+1}(B\rtimes_r G)\]
such that
\begin{enumerate}
\item[(i)] $J_{red,G}(z)$ induces right multiplication by $j_{red,G}(z)$ in $K$-theory ;
\item[(ii)] $J_{red,G}$ is additive, i.e.
\[J_{red,G}(z+z')=J_{red,G}(z)+J_{red,G}(z').\]
\item[(iii)] For every $G$-morphism $f : A_1\rightarrow A_2$,
\[J_{red,G}(f^*(z))=J_{red,G}(z)\circ f_{G,red,*}\] for all $z\in KK_1^G(A_2,B)$.
\item[(iv)] For every $G$-morphism $g : B_1\rightarrow B_2$,
\[J_{red,G}(g_*(z))= g_{G,red,*}\circ J_{red,G}(z)\] for all $z\in KK_1^G(A,B_1)$.
\item[(v)] Let $0\rightarrow J\rightarrow A\rightarrow A/J\rightarrow 0$ be a semi-split equivariant extension of $G$-algebras and $[\partial_J]\in KK_1^G(A/J,J)$ be its boundary element. Then 
\[J_G([\partial_J])=D_{J\rtimes_r G,A\rtimes_rG}.\] 
\end{enumerate}
\end{prop}

\begin{dem}
\begin{enumerate}

\item[(i)]The $K$-cycle $[\partial_{K_{B\rtimes_r G},E^{(\pi,T)}}]\in KK_1(A\rtimes_r G, B\rtimes_r G)$ implementing the boundary of the extension $E^{(\pi,T)}$ induces the map $-\otimes_{A\rtimes_r G} j_{red,G}(z)$ by definition, and modulo Morita equivalence, which immediately gives the first point.

\item[(ii)] If $z,z'$ are elements of $KK_1^G(A,B)$, represented by two $K$-cycles $(H_B,\pi_j,T_j)$, let $(\tilde H_B,\pi,T)$ be $(H_B\oplus H_B,\pi_0\oplus \pi_1,T_0\oplus T_1)$ which is a $K$-cycle representing the sum $z+z'$. Then $E^{(\pi,T)}$ is naturally isomorphic to the extension sum of the $E_j:=E^{(\pi_j,T_j)}$, namely
\[\begin{tikzcd}[column sep = small]
0\arrow{r} & \mathfrak M_2(K_{B\rtimes_r G}) \arrow{r} & D \arrow{r} & A\rtimes_r G \arrow{r} & 0
\end{tikzcd}\]
where 
\[D=\left\{\begin{pmatrix}x_1 & k_{12}\\ k_{21} & x_2\end{pmatrix} : x_j\in E_j , p_1(x_1)=p_2(x_2), k_{ij}\in K_{B\rtimes_r G}\right\},\]
with $p_j : E_j\rightarrow A\rtimes_r G$ the $*$-homomorphisms of the extensions.
Naturality of the controlled boundary maps (\cite{OY2}, remark $3.5.(i)$) ensures that the boundary of the sum of two extensions is the sum of the boundary of each, thus the result holds.
\item[(iii)] Let $z\in KK_1^G(A_2,B)$, represented by a cycle $(H_B,\pi,T)$. Representing $f^*(z)$ is $(H_B,f^*\pi,T)$ with $f^*\pi=\pi \circ f$. The map 
\[\phi : \left\{\begin{array}{lll} E^{f^*(\pi,T)} & \rightarrow & E^{(\pi,T)} \\
( x, P_G(f^*\pi)(x)P_G+y) & \rightarrow & ( f_G(x), P_G(f^*\pi)(x)P_G+y) \end{array}\right. \]
satisfies
\begin{enumerate}
\item[$\bullet$] $\phi(K_{B\rtimes_r G})\subseteq K_{B\rtimes_r G}$, and makes the following diagram commute
\[\begin{tikzcd}[column sep = small]
0\arrow{r} & K_{B\rtimes_r G}\arrow{r}\arrow{d}{=} & E^{f^*(\pi,T)} \arrow{r}\arrow{d}{\phi}& A_1\rtimes_r G\arrow{r}\arrow{d}{f_G} & 0\\
0\arrow{r} & K_{B\rtimes_r G}\arrow{r} & E^{(\pi,T)} \arrow{r}& A_2\rtimes_r G\arrow{r} & 0
\end{tikzcd}.\]
\item[$\bullet$] It intertwines the sections of the two extensions.
\end{enumerate}
Remark \ref{rk3.8} ensures that \[D_{K_{B\rtimes_r G}, E^{f^*(\pi,T)} } =  D_{K_{B\rtimes_r G}, E^{(\pi,T)} }\circ f_{G,*},\] and the claim is clear from composition by $\mathcal M_{B\rtimes_r G}^{-1}$.

\item[(iv)] %New proof.\\
Let $z \in KK^G(A,B_1)$ be represented by the $K$-cycle $(H_{B_1},\pi,T)$. Let $V\in \mathcal L_{B_2}(H_{B_1}\otimes_g B_2,H_{B_2})$ be the isometry of remark \ref{isometry}. Notice that $V$ intertwines the actions of $G$ on $ H_{B_1}\otimes B_2 $ and $H_{B_2}$. According to Lemma \ref{isometryKK}, 
\[g_*(z)=[H_{B_1}\otimes_g B_2, \pi\otimes_g 1, T\otimes_g 1]\in KK^G(A,B_2)\] 
is also represented by $[H_{B_2}, \pi',T' ]$ where $\pi' = Ad_{V}\circ (\pi\otimes_g 1)$ and $T' = V(T\otimes_g 1)V^* +1-VV^*$. Let $\psi$ be given by the composition $Ad_{V_G}\circ g_G$, and $P' = \frac{1+T'}{2}$.\\
The map $\Psi :(x,y)\mapsto (x, \psi(y))$ defines a $*$-homomorphism $E^{(\pi,T)} \rightarrow E^{(\pi',T')}$ such that 
\[\Psi(x,P_G\pi_G(x)P_G +y)= (x, P'_G  \pi_G'(a)P'_G + \psi(y)) ).\] 
Indeed, the crossed-product functor commutes with pull-back by $G$-morphisms, and $Ad_{V_G}\circ g_G \circ\pi_G= (Ad_V\circ g_* \circ \pi)_G = \pi'_G$ and $\psi(P_G)= V_G (P_G\otimes_{g_G} 1)V^*_G = (V(P\otimes_g 1 ) V^*)_G = (P')_G$ so that 
\[\psi(P_G \pi_G(x) P_G)=P'_G \pi'_G(x) P'_G. \]
This gives a commutative diagram 
\[\begin{tikzcd}[column sep = small]
0\arrow{r} & K_{B_1\rtimes G}\arrow{r}\arrow{d}{\psi} & E^{(\pi,T)} \arrow{r}\arrow{d}{\Psi}& A\rtimes_r G\arrow{r}\arrow{d}{=} & 0\\
0\arrow{r} & K_{B_2\rtimes G}\arrow{r} & E^{(\pi',T')} \arrow{r}& A\rtimes G\arrow{r} & 0
\end{tikzcd}.\]
and $\Psi$ intertwines the two filtered sections by the previous relation. Moreover, recall from remark \ref{isometry} that $\Psi( K_{B_1\rtimes G})\subseteq K_{B_2\rtimes G}$, hence we can again apply the remark \ref{rk3.8} to state
\[ D_{K_{B_2\rtimes G},E^{(\pi',T')}}=\psi_*\circ D_{K_{B_1\rtimes G},E^{(\pi,T)}},\]
which we compose by the Morita equivalence on the left $M_{B_2\rtimes G}^{-1}$
\[J_G(g_*(z)) = M_{B_2\rtimes G}^{-1}\circ g_{G,*}\circ D_{K_{B_1\rtimes G},E^{(\pi,T)}}.\]
Notice that $\psi$ is the $*$-homomorphism induced by $g_G$ on the compact operators, hence the homomorphisms inducing the Morita equivalence make the following diagram commutes,
\[\begin{tikzcd}B_1\rtimes G\arrow{r}{g_G}\arrow{d} & B_2\rtimes G\arrow{d} \\ K_{B_1\rtimes G } \arrow{r}{\psi}& K_{B_2\rtimes G }\end{tikzcd},\]
and $J_G(g_*(z))= g_{G,*}\circ M_{B_1\rtimes G}^{-1}\circ D_{K_{B_1\rtimes G},E^{(\pi,T)}}=g_{G,*}\circ J_G(z)$.\\

%% NEW NEW PROOF
\item[(v)] We can suppose $A$ unital. Let $0 \rightarrow J \rightarrow A \rightarrow A /J \rightarrow 0$ be a $G$-equivariant semi-split extension of $G$-algebras with $q:A\rightarrow A/J$ the quotient map.. Let us denote by $s : A/J \rightarrow A $ the $G$-equivariant completely positive cross section.   According to the equivariant Kasparov-Stinespring theorem, %\ref{GKasparovStinespring}
there exists an equivariant $A$-Hilbert module $E$ and a $G$-equivariant $*$-homomorphism $\pi : A/J \rightarrow \mathcal L_{A}(A\oplus E)$ such that $s(x) = P \pi(x) P$, where $P \in \mathcal L_{A}(A\oplus E)$ is the projection on the $A$ factor. Consider the $J$-Hilbert module $E' = (A\oplus E)\otimes_J J \cong J\oplus (E\otimes_J J)$, and the natural map $\tilde\pi =\pi\otimes_J 1: A/ J \rightarrow \mathcal L_{J}(E')$. Put $\tilde T= (2P-1)\otimes_J id_J$. By the stabilization theorem, we can suppose that $E'$ is a standard $G$-equivariant $J$-Hilbert module. 
Put 
\[\psi_0(x)(y \oplus \xi ) =  (xy) \oplus 0 \quad \forall \xi \in E\otimes J,\forall y\in J, \] 
for every $x\in A$. This defines a $G$-equivariant $*$-homomorphism $\psi_0 : A \rightarrow \mathcal L_J(E')$ such that $\psi_0(x)\in  \mathfrak K _J (E')$ when $x\in J$. Put $\psi ( a ) = (q(a), \psi_0 (a))$. As $P\tilde \pi (q(a))P -\psi_0(a) = \psi_0( s(q(a)) - a )\in \mathfrak K_J(E')$, $\psi(a)\in E^{(\tilde\pi,\tilde T)}$ holds, and the following diagram is commutative with exact rows:
\[\begin{tikzcd}[column sep = small]
0\arrow{r} &J\arrow{r}\arrow{d}{\psi_0} & A \arrow{r}\arrow{d}{\psi}& A/J\arrow{r}\arrow{d}{=} & 0\\
0\arrow{r} & \mathfrak K_J(E')\arrow{r} & E^{(\tilde \pi, \tilde T)} \arrow{r}& A/J\arrow{r} & 0
\end{tikzcd},\]
Hence, by functoriality and semi-split exactness of the reduced cross product, the following diagram commutes
\[\begin{tikzcd}[column sep = small]
0\arrow{r} & J\rtimes_r G\arrow{r}\arrow{d}{(\psi_0)_G} & A\rtimes_r G \arrow{r}\arrow{d}{\psi_G}& A/J\rtimes_r G \arrow{r}\arrow{d}{=} & 0\\
0\arrow{r} & \mathfrak K \otimes J\rtimes_r G\arrow{r} & E^{(\tilde \pi, \tilde T)} \arrow{r}& A/J\rtimes_r G \arrow{r} & 0
\end{tikzcd},\]
and remark \ref{rk3.8} ensures that $J_G([\partial_J]) = D_{J\rtimes_r G,A\rtimes_r G}$. 
%% END NEW NEW PROOF

\qed
\end{enumerate}
\end{dem}

\subsubsection{Even case}

We can now define $J_G$ for even $K$-cycles. Let $A$ and $B$ be two $G$-algebras. Let $[\partial_{SB}]\in KK_1(B,SB)$ be the $K$-cycle implementing the boundary of the extension $0\rightarrow SB\rightarrow CB\rightarrow B\rightarrow 0$, and $[\partial]\in KK_1(\C,S)$ be the Bott generator. As $z\otimes_B [\partial_{SB}]$ is an odd $K$-cycle, we can define
\[J_G(z):= \hat\tau_{B\rtimes G}([\partial]^{-1})\circ J_G(z\otimes[\partial_{SB}]).\] 

Here $\hat\tau_D(z)$ refers, for any $C^*$-algebras $D,A_1,A_2$ and $z\in KK_*(A_1,A_2)$, to the $(\alpha_\tau,k_\tau)$-controlled map $\hat K (A_1\otimes D )\rightarrow \hat K(A_2\otimes D)$ of theorem \ref{tensorization}. We can see that, if we set $\alpha_J=\alpha_\tau \alpha_D$ and $k_J=k_\tau * k_D$, $J_G(z)$ is $(\alpha_J,k_J)$-controlled.\\

\begin{prop}\label{Kasparov}
Let $A$ and $B$ two $G$-$C^*$-algebras. For every $z\in KK^G_*(A,B)$, there exists a control pair $(\alpha_J,k_J)$ and a $(\alpha_J,k_J)$-controlled morphism
\[J_{red,G}(z) : \hat K(A\rtimes_r G)\rightarrow \hat K(B\rtimes_r G)\]
of the same degree as $z$, such that
\begin{enumerate}
\item[(i)] $J_{red,G}(z)$ induces right multiplication by $j_{red,G}(z)$ in $K$-theory ;
\item[(ii)] $J_{red,G}$ is additive, i.e.
\[J_{red,G}(z+z')=J_{red,G}(z)+J_{red,G}(z').\]
\item[(iii)] For every $G$-morphism $f : A_1\rightarrow A_2$,
\[J_{red,G}(f^*(z))=J_{red,G}(z)\circ f_{G,red,*}\] for all $z\in KK_*^G(A_2,B)$.
\item[(iv)] For every $G$-morphism $g : B_1\rightarrow B_2$,
\[J_{red,G}(g_*(z))= g_{G,red,*}\circ J_{red,G}(z)\] for all $z\in KK_*^G(A,B_1)$.
\item[(v)] $J_{red,G}([id_A]) \sim_{(\alpha_J,k_J)} id_{\hat K(A\rtimes G)}$
\end{enumerate}
\end{prop}

\begin{dem} Let $z\in KK_0^G(A,B)$.
\begin{itemize}
\item[$(i)$] $J_{red,G}(z)$ induces in $K$-theory right multiplication by $j_{red,G}(z\otimes [\partial_{SB}])\otimes \tau_{B\rtimes_r G}([\partial]^{-1})$. But $j_{red, G}$ respects Kasparov products, and by \cite{LeGall}, $j_{red,G}([\partial_{SB}]) = [\partial_{SB\rtimes_r G}]$ and $\tau_{B\rtimes_r G}([\partial]^{-1}) = [\partial_{SB\rtimes_r G }]^{-1}$ are inverse of each others.
\item[$(ii)$] Additivity follows from additivity of the Kasparov product and of proposition \ref{Kasparov1}.
\item[$(iii)$] Recall the equality $f^*(x)\otimes y = f^*(x\otimes y)$. As a consequence of the previous proposition \ref{Kasparov1},
\[\begin{array}{lcl} J_{red,G}(f^*(z)) & = &  \hat\tau_{B\rtimes_r G}([\partial]^{-1}) \circ J_{red,G}(f^*(z)\otimes [\partial_{SB}]) \\
		& = &  \hat\tau_{B\rtimes_r G}([\partial]^{-1}) \circ J_{red,G}(f^*(z\otimes [\partial_{SB}]))\\
		& = &  \hat\tau_{B\rtimes_r G}([\partial]^{-1}) \circ J_{red,G}(z\otimes [\partial_{SB}])\circ f_{G,red,*}\\
		& = & J_{red,G}(z) \circ f_{G,red,*}
\end{array}\]
\item[$(iv)$] We have $g_*(z)\otimes_{B'} [\partial_{SB'}] = z\otimes_B g^{*}[\partial _{SB'}] = (Sg)_{*}(z\otimes_B [\partial _{SB}])$, hence, by proposition \ref{Kasparov1}, 
\[J_{red,G}(g_*(z)) = \hat\tau_{B'\rtimes_r G}([\partial]^{-1}) \circ (Sg)_{G,red,*} \circ J_{red,G}(z\otimes_B [\partial_{SB}]),\]
but, using properties of $\hat\tau_{SB}$, see \ref{tensorization}, $\tau_{B'\rtimes_r G}([\partial]^{-1})\circ (Sg)_{*,red,G} = g_{G,red,*}\circ\tau_{B\rtimes_r G}(([\partial]^{-1})$.
\item[$(v)$] By definition, $J_G([id_A]) = \hat\tau_{A\rtimes_r G}([\partial]^{-1}]) \circ J_G([\partial_{SA}])$. By point $(v)$ of \ref{Kasparov1}, $J_G([\partial_{SA}]) = D_{A\rtimes_r G,SA\rtimes_r G}$ and by property of $\hat \tau$, $\hat\tau_{A\rtimes_r G}([\partial]^{-1}]) $ is a controlled left inverse of $\hat\tau_{A\rtimes_r G}([\partial]])= D_{A\rtimes_r G,SA\rtimes_r G}$, hence $J_G([id_A])\sim id_{\hat K(A\rtimes_r G)}$. 
\end{itemize}
\qed
\end{dem}

We now show that the controlled Kasparov transform respects in a quantitative way the Kasparov product.

\begin{prop} There exists a control pair $(\alpha_K,k_K)$ such that for every $G$-$C^*$-algebras $A$, $B$ and $C$, and every $z\in KK^G(A,B),z'\in KK^G(B,C)$, the controlled equality
\[J_G(z\otimes_B z') \sim_{\alpha_K,k_K} J_G(z')\circ J_G(z)\]
holds.
\end{prop}
\begin{dem}

Recall from \ref{propertyD} that every $\alpha\in KK_0^G(A,B)$ satisfies property $(d)$ for a universal $d$. By naturality, this property reduces the proof to the special case of $\alpha\in KK^G(A,B)$ being the inverse of a $*$-homomorphism $\theta : B\rightarrow A$ in $KK^G$-theory : $\alpha\otimes_B [\theta]=1_A$. Let $z\in KK^G(B,C)$ :
\[\begin{array}{rcl}
J_G (\alpha\otimes z) & \sim_{\alpha_J^2,k_J*k_J} &  J_G(\alpha\otimes z)\circ J_G(\alpha\otimes [\theta]) \\
			& \sim & J_G(\alpha\otimes z)\circ J_G(\theta_*(\alpha))\\
			& \sim & J_G(\alpha\otimes z)\circ \theta_{G,*}\circ J_G(\alpha)\\
			& \sim & J_G(\theta^*(\alpha\otimes z))\circ J_G(\alpha)\\
			& \sim & J_G(z)\circ J_G(\alpha) \\
\end{array}\] 
because $\theta^*(\alpha\otimes z)=\theta^*(\alpha)\otimes z=1\otimes z =z$. The control on the propagation of the first line follows from remark $2.5$ of \cite{OY2} and point $(v)$, the other lines are equal by points $(iii)$ and $(iv)$. As $d$ is uniform for all locally compact groupoids with Haar systems, a simple induction concludes, and $(\alpha_K,k_K)$ can be taken to be $( \alpha_J^{2d},( k_J)^{*2d})$. The same argument works if $z'$ is an even $KK$-element.\\

Let $z$ and $z'$ be odd $KK$-elements. Then we have :
\[\begin{array}{rcl}
J_G (z\otimes z') & = &  J_G (z\otimes_B [\partial_{B}]\otimes_{SB} [\partial_B]^{-1}\otimes_B z') \\
			& \sim & J_G ( [\partial_B]^{-1}\otimes_B z')\circ J_G (z\otimes_B [\partial_{B}])\\
			& \sim & J_G ( [\partial_B]^{-1}\otimes_B z')\circ J_G ( [\partial_B])\circ J_G( [\partial_B]^{-1})\circ J_G (z\otimes_B [\partial_{B}])\\		
			& \sim & J_G (  z')\circ J_G (z),\\
\end{array}\] 
where we used the previous case for the second line, Lemma \ref{Kasparov1} for the third line, and Proposition \ref{InverseEven} for the last one.\\
\qed
\end{dem}

%%%%%%%%%%%%%%%%%%%%%%%%%%%%%%%%%%%%%%%
%%%%%%%%%%%%%%%%%%%%%%%%%%%%%%%%%%%%%%%
\subsection{Controlled assembly maps}
%%%%%%%%%%%%%%%%%%%%%%%%%%%%%%%%%%%%%%%
%%%%%%%%%%%%%%%%%%%%%%%%%%%%%%%%%%%%%%%

Let $E\in\mathcal E$. Recall that every $\eta\in P_E(G)$ is a finite probability measure in $G^{p(\eta)}$, hence can be written as a finite sum 
\[\eta = \sum_{g\in G^{p(\eta)}}\lambda_g(\eta)\delta_g,\]
where $\delta_g$ is the Dirac probability at $g\in G^{p(\eta)}$. This defines a map 
\[\lambda : \left\{\begin{array}{rcl}
G\times_{s,p} P_E(G) & \rightarrow & [0,1] \\
 (g,\eta) & \mapsto & \lambda_g(\eta)\end{array}\right. .\] 

\begin{lem}
The function $\lambda : G\times_{s,p} P_E(G) \rightarrow [0,1]$ is continuous. 
\end{lem}

\begin{dem}
Let $\{(g_j,\eta_j)\}_{j\in J}\subseteq G\times_{s,p}P_E(G)$ be a generalized sequence converging to $(g,\eta)\in G\times_{s,p}P_E(G)$. Let $x=p(\eta) = s(g)$ and $x_j=p(\eta_j) = s(g_j)$. There exists an open neighborhood $U\subseteq G^{(0)}$ of $x$ such that we can decompose $s^{-1}(U) = \coprod_{g\in G^x} S_g$ with, for each $g\in G^x$, an open bisection $S_g\subseteq G$ and such that $x_j$ belongs ultimately to $U$. Notice that $s$ is injective on each $S_g$. We can suppose $x_j\in U$ for every $j$. Let, for every $g\in G^x$, $\phi_g : G\rightarrow [0,1]$ be a continuous function equal to $1$ on $S_g$ and vanishing outside a small neighborhood of $S_g$. Then $\lambda_g(\eta) =\langle \eta,\phi_g\rangle$ and $\lambda_{g_j}(\eta_j) = \langle \eta_j,\phi_g\rangle$ for every $j\in J$. It is now clear, by definition of the weak-$*$ topology, that $\lambda_{g_j}(\eta_j)$ converges to $\lambda_{g}(\eta)$.  \\
\qed %$\langle \eta,\phi_g\rangle / \sum_{h\in G^x}\langle \eta,\phi_h\rangle$ / \sum_{h\in G^x}\langle \eta_j,\phi_h\rangle
\end{dem}

%Let $K\subseteq P_E(G)$ be a compact. Then $\{g\in G \text{ s.t. }\lambda(g,\eta)= 0\ \forall \eta \in K\}$ is compact. 
Notice that $\lambda_{g'}(g.\eta) = \lambda_{g^{-1}g'}(\eta)$ for every $\eta\in P_E(G)$, $g\in G_{p(\eta)}$ and $g'\in G^{r(g)}$. Define $h$ as the restriction of $\lambda$ to $G^{(0)}\times_{s,p}P_E(G)$, i.e. $h(x,\eta) = \lambda_{e_x}(\eta)$, for every $x\in G^{(0)}$ and $\eta\in P_E(G)$ such that $p(\eta) = x$. Notice that $h(x,g^{-1}.\eta)= \lambda(g,\eta) = \lambda_g(\eta)$ for every $(g,\eta)\in G\times_{s,p} P_E(G)$. Let us set
 
\[\mathcal L_E(g,\eta) = h^{\frac{1}{2}}(x,\eta) h^{\frac{1}{2}}(x,g^{-1}.\eta)= \lambda_{e_x}^{\frac{1}{2}}(\eta)\lambda_g^{\frac{1}{2}}(\eta),\] 

where $x= r(g)$. Notice that $\text{supp }\mathcal L_E\subseteq E\times_{s,p} p^{-1}(r(E))$, which is compact. Hence $\mathcal L_E$ is an element of $C_c(G,C_0(P_E(G)))$. Because $\sum_{g\in G^{p(\eta)}}\lambda_g(\eta) = 1 ,\forall (g,\eta)\in G\times_{s,p} P_E(G)$, this defines a projection in $C_0(P_E(G))\rtimes_r G$ of propagation less than $E$. Moreover $\mathcal L_E$ defines a controlled $K$-theory class $[\mathcal L_E,0]_{\varepsilon,F}$ for any $\varepsilon\in (0,\frac{1}{4})$ and any $F\in \mathcal E$ such that $E\subseteq F$.\\ %Moreover, for such controlled subsets, if $q_E^F : C_0(P_E(G))\rightarrow C_0(P_F(G))$, then $(q_E^F)_*[\mathcal L_E,0]_{\varepsilon,E} = [\mathcal L_F,0]_{\varepsilon,F}$.\\

%There exists continuous functions $\lambda_g : P_E(G)\rightarrow [0,1]$ for all $g\in G$ such that $\eta = \sum_{g\in G^x}\lambda_g(\eta) \delta_g$ for any $\eta\in P_E(G)$. Define $h(x)=\lambda_{e_x}$ and $\phi = \sqrt h$. Notice that $g.h = \lambda_g$ for all $g\in G$, and as $\sum_{g\in G^x}\lambda_g = 1 ,\forall x\in G^{(0)}$, 
%\[\mathcal L_E =\sum_{g\in G^x} \phi(g.\phi)\]
%defines a projection of $C_0(P_E(G))\rtimes_r G$ with bounded propagation, and defines a $K$-theory class $[\mathcal L_E,0]_{\varepsilon,F}$ for any $\varepsilon\in (0,\frac{1}{4})$ and any $F\in \mathcal E$ such that $E\subseteq F$. Moreover, for such controlled subsets, if $q_E^F : C_0(P_E(G))\rightarrow C_0(P_F(G))$, then $(q_E^F)_*[\mathcal L_E,0]_{\varepsilon,E} = [\mathcal L_F,0]_{\varepsilon,F}$.\\

For every $G$-algebra $B$ and every controlled subsets $E,E'\in\mathcal E$ such that $E\subseteq E'$, the canonical inclusion $P_E(G)\hookrightarrow P_{E'}(G)$ induces a $*$-homomorphism $q_E^{E'} : C_0(P_{E'}(G))\rightarrow C_0(P_{E}(G))$, hence a map $(q_E^{E'})^* : RK^G(P_E(G),B)\rightarrow RK^G(P_{E'}(G),B)$ in $K$-homology and a map $((q_E^{E'})_G)_* : K(C_0(P_{E'}(G))\rtimes_r G)\rightarrow K(C_0(P_{E}(G))\rtimes_r G)$ in $K$-theory. The family of projections $\mathcal L_E$ are compatible with the morphisms $q_E^{E'}$, i.e. $((q_E^{E'})_G)_*[\mathcal L_{E'},0]_{\varepsilon,E'} = [\mathcal L_{E},0]_{\varepsilon,E}$, for every $\varepsilon\in (0,\frac{1}{4})$.

\begin{definition}
Let $B$ be a $G$-algebra, $\varepsilon\in (0,\frac{1}{4})$, and $E\in\mathcal E$. Let $F\in \mathcal E$ such that $k_J(\varepsilon).E \subseteq F$. The controlled assembly map for $G$ is defined as the family of maps :
\[\mu_{G,B}^{\varepsilon,E,F}\left\{
\begin{array}{rcl}
RK^G(P_E(G), B) & \rightarrow & K_*^{\varepsilon, F}(B\rtimes_r G)\\
z & \mapsto & \iota_{\alpha_J\varepsilon', k_J(\varepsilon').F'}^{\varepsilon,F} \circ J_G^{\varepsilon', F'}(z)([\mathcal L_E,0]_{\varepsilon' , F'})
\end{array}\right.\]
where $\varepsilon'$ and $F'$ satisfy :
\begin{itemize}
\item[$\bullet$] $\varepsilon'\in (0,\frac{1}{4})$ such that $\alpha_J \varepsilon'\leq \varepsilon$,
\item[$\bullet$] and $F'\in\mathcal E$ such that $E\subseteq F'$ and $k_J(\varepsilon').F'\subseteq F$.
\end{itemize}
\end{definition}

%%%%%%%%%%%%
%%%%%%%%%%%%
%%%%%%%%%%%%

\begin{rk}
The assembly map is defined for any reasonnable crossed-products by $G$. In particular for the reduced one and the maximal one, so that we have two different assembly maps, which we shall distinguish writing $\mu_{G,r}$ and $\mu_{G,max}$ if necessary.
\end{rk}

\begin{rk} As for the controlled coarse assembly map, the controlled assembly map for an étale groupoid is compatible with the structure morphisms $q_{E}^{E'}$, and $\iota_{\varepsilon,E}^{\varepsilon',F'}$.
\end{rk}

\begin{rk} 
The family of assembly maps $\mu_{G,B}^{\varepsilon,E,F}$ induces the Baum-Connes assembly map for $G$ in $K$-theory. Notice that $h$ is a cutoff function for the action on $P_E(G)$, hence $\mathcal L_E$ coincides with $\mathcal L_{P_E(G)}$. Moreover, the following diagram commutes
\[\begin{tikzcd}
RK_*^G(P_E(G),B) \arrow{r}{\mu_{G,B}^{\varepsilon,E,F}}\arrow{dr}{\mu_{G,B}^E} & K_*^{\varepsilon, F}(B\rtimes_r G)\arrow{d}{\iota_{\varepsilon,F}}\\ 
		&  K_*(B\rtimes_r G)
\end{tikzcd}\]
because $J_G(z)$ induces the right multiplication by $j_G(z)$ and also $\mu_G^E(z)=[\mathcal L_E]\otimes j_G(z)$ by remark \ref{projection}. But, as $(q_E^{E'})_*[\mathcal L_E,0]_{\varepsilon,E} = [\mathcal L_{E'},0]_{\varepsilon,E'}$ as soon as $E\subseteq E'$, this diagram commutes with inductive limit over $E$.\\
\end{rk}

%%%%%%%%%%%%%%%%%%%%%%%%%%%%%%%%%%%%
\subsection{Controlled Baum-Connes conjecture}
%%%%%%%%%%%%%%%%%%%%%%%%%%%%%%%%%%%%

We generalize in this subsection to the setting of étale groupoids the results obtained by H. Oyono-Oyono and G. Yu in \cite{OY2}. We will define a controlled version of the Baum-Connes conjecture, and give conditions that ensure it is satisfied. In particular, when $G$ satisfies the Baum-Connes conjecture with coefficients, it satisfies the controlled Baum-Connes conjecture. To this end, we will have to compute the $K$-homology of a finite typed $G$-simplicial complex with coefficients in the infinite product of stable $G$-algebras.\\

In this subsection, the groupoid $G$ must have a compact unit space $G^{(0)}$. Then, if $\{B_j\}$ is a countable family of $G$-algebras, the diagonal action of $G$ on the product provides $\prod_j B_j$ with the structure of a $G$-algebra. \\

Let $A$ be a $G$-algebra. We will say that :\\

$\bullet$(Quantitative Injectivity) $\mu_{G,A}$ is quantitatively injective if, for every $E\in \mathcal E$, there exists $\varepsilon\in (0,\frac{1}{4})$ such that, for every $F\in\mathcal E$ satisfying $k_J(\varepsilon). E \subseteq F$, there exists $E'\in\mathcal E$ such that $E\subseteq E'$ and 
\[\forall z\in RK^G(P_E(G),A),\quad\mu_{G,A}^{\varepsilon,E,F}(z) = 0 \text{ implies } q_E^{E'}(z)=0.\]

$\bullet$(Quantitative Surjectivity) $\mu_{G,A}$ is quantitatively surjective if there exists $\varepsilon$ such that, for every $F\in \mathcal E$ such that, there exists $\varepsilon' \in (\varepsilon,\frac{1}{4})$ and $E,F'\in\mathcal E$ such that $F \subseteq F'$ and $k_J(\varepsilon'). E\subseteq F'$,  
\[\forall y\in K^{\varepsilon, F}(A\rtimes_r G),\exists z\in RK^G(P_E(G),A) \text{ such that } \mu_{G,A}^{\varepsilon', E, F'}(z)=\iota_{\varepsilon,F}^{\varepsilon',F'}(y).\]

\begin{prop} 
Let $G$ be an étale groupoid with compact base space, and let $B$ be a $G$-algebra.\\

$\bullet$(Quantitative Injectivity) If $\mu_{G,A}$ is quantitatively injective then $\mu_{G,A}$ is one-to-one.\\
%$\bullet$(Quantitative Surjectivity) there exists $\varepsilon\in (0,\frac{1}{4})$ such that for any controlled subsets $E,F\in\mathcal E$  such that $k_J(\varepsilon).E\subseteq F$,$\exists \varepsilon',F'$ such that $\varepsilon'\leq \varepsilon<\frac{1}{4}$ and $k_J E\subseteq F\subseteq F'$, such that for all $y\in K_*^{\varepsilon,F}(A\rtimes_r G),\exists x \in RK_*^G(P_E(G),A)$ such that $\mu_{G,A}^{\varepsilon',E,F'}(x)=\iota_{\varepsilon,F}^{\varepsilon',F'}(y)$;\\
$\bullet$(Quantitative Surjectivity) If $\mu_{G,A}$ is quantitatively surjective then $\mu_{G,A}$ is onto.\\
\end{prop}

\begin{dem}
Let $E\in\mathcal E$ and $x\in RK(P_E(G)),A)$ such that $\mu_{G,A}^E(x)=0$. Then, for every $\varepsilon\in (0,\frac{1}{4})$ and every $F\in\mathcal E$ such that $k_J(\varepsilon).E\subseteq F$, $\iota_{\varepsilon,F}\circ\mu_{G,A}^{\varepsilon,E,F}(x)=0$. Let $\varepsilon''>0$ and $F''\in\mathcal E$ satisfying $\alpha_J\varepsilon''\leq \varepsilon$ and $k_J(\varepsilon''). F''\subseteq F$. By remark \ref{approximation}, there exists a universal $\lambda\geq 1$ and a controlled subset $F'\in\mathcal E$ such that $F\subseteq F'$ and
\[\begin{array}{lll}0 &  =  & \iota_{\varepsilon,F}^{\lambda\varepsilon,F'}\circ \mu_{G,A}^{\varepsilon,E,F}(x) \\
			& = & \iota_{\varepsilon,F}^{\lambda\varepsilon,F'}\circ \iota^{\varepsilon,F}_{\alpha_J\varepsilon'',k_J(\varepsilon'').F''} (J_{G}^{\varepsilon'',F''}(x)([\mathcal L_E,0]_{\varepsilon'',F''})) \\
			& = & J_{G}^{\lambda\varepsilon,F'}(x)([\mathcal L_E,0]_{\lambda\varepsilon,F'}) \\
			& = & \mu_{G,A}^{\lambda\varepsilon,E,F'}(x).
\end{array}\]
But then the quantitative injectivity condition ensures that $q_E^{E'}(x)=0$ in $RK^G(P_{E'}(G),A)$ and $x=0$ in $K^{top}(G,A)$, which is an inductive limit over $E$.\\

Let us prove the second point. Let $y\in K(B\rtimes_r G)$, and let $\varepsilon\in (0,\frac{1}{4})$, $F\in\mathcal E$ and $x\in K^{\varepsilon,F}(B\rtimes_r G)$ such that $\iota_{\varepsilon,F}(x) =y$. The quantitative surjectivity condition implies that there exist $\varepsilon'\in (0,\frac{1}{4})$, $E,F'\in\mathcal E$, and $z\in RK^G(P_E(G),B)$ satisfying $\varepsilon\leq \varepsilon'$, $k_J(\varepsilon').E\subseteq F'$, $F\subseteq F'$ and $\mu_{G,B}^{\varepsilon',E,F'}(z) = \iota_{\varepsilon,F}^{\varepsilon',F'} (x)$, hence $\mu^E_{G,B}(z) =y$.\\
\qed
\end{dem}

This kind of statement leads us to define the following properties, following \cite{OY3}.\\
\begin{itemize}
\item[$\bullet$] $QI_{G,B}(E,E',F,\varepsilon)$ : for any $x\in RK^G(P_E(G), B )$, then $\mu^{\varepsilon,E,F}_{G,B}(x) = 0$ implies $q_E^{E'}(x)=0$ in $RK^G(P_{E'}(G),B)$.
\item[$\bullet$] $QS_{G,B}(E,F,F',\varepsilon,\varepsilon')$ : for any $y\in K^{\varepsilon,F}(B\rtimes_r G)$, there exists $x\in RK^G(P_E(G),B)$ such that $\mu^{\varepsilon',E,F'}_{G,B}(x)=\iota_{\varepsilon,F}^{\varepsilon',F'}(y)$.\\
\end{itemize} 

Controlled injectivity and surjectivity being defined, we can now state a controlled version of the Baum-Connes conjecture.

\begin{definition} 
Let $\lambda \geq 1$. The groupoid $G$ is said to satisfy the controlled Baum-Connes conjecture, with coefficients in the $G$-algebra $B$, at scale $\lambda$ if:
\begin{itemize}
\item[$\bullet$] for every $\varepsilon \in (0,\frac{1}{4\lambda})$ and $E,F\in \mathcal E$ such that $k_J(\varepsilon).E \subseteq F$, there exists $E'\in \mathcal E$ such that $E\subseteq E'$ and $QI_{G,B}(\varepsilon, E,E',F)$ holds;
\item[$\bullet$] for every $\varepsilon \in (0,\frac{1}{4\lambda})$ and $F\in \mathcal E$, there exist $E,F'\in \mathcal E$ such that $F \subseteq F'$, $k_J(\varepsilon).E \subseteq F'$ and $QS_{G,B}(\varepsilon,\lambda \varepsilon , E,F,F')$ holds. 
\end{itemize}
\end{definition}

The main result of this section is the following theorem.

\begin{thm} \label{ControlledBC} There exists a universal $\lambda \geq 1$ such that for every étale groupoid with compact base space $G$, and every $B$ $G$-algebra, the following statements are equivalent:
\begin{itemize} 
\item[$\bullet$] $G$ satisfies the controlled Baum-Connes conjecture with coefficients in $B$ at scale $\lambda$, 
\item[$\bullet$] $G$ satisfies the Baum-Connes conjecture with coefficients in $l^\infty (\N,B\otimes \mathfrak K)$.
\end{itemize}
\end{thm}

%%%%%%%%%%%%%%%%%%%%%%%%%%%%%%%%%%%%%%%%%%%%%%%%%%%%%%%%%%%%%%%%%%%%%%%%%%%%%
\subsection{Products in topological $K$-theory and in controlled $K$-theory}
%%%%%%%%%%%%%%%%%%%%%%%%%%%%%%%%%%%%%%%%%%%%%%%%%%%%%%%%%%%%%%%%%%%%%%%%%%%%%

To prove the main theorems of this section, what we called quantitative statements, we will need a serie of lemmas. The first part is devoted to a study of the behaviour of $K$-homology under product of the coefficients. The strategy is very similar to the proofs of the part $3$ of \cite{TuBC2}. \\

%Let us begin with a remark that will be useful later.

%\begin{rk}\label{construction}
%Let $\{B_j\}_{j\in\N}$ be a family of $G$-algebras and let $T_j \in \mathfrak K (H_{B_j})\cong B_j \otimes\mathfrak K$ such that $\sup_j ||T_j||<\infty$. Then we can define  
%\end{rk}

Recall first that the Haar system $(\lambda^x)_{x\in G^{(0)}}$ on $G$ defines, for every $f\in C_c(G)$, a continuous function on the base space $x\mapsto \int_{g\in G^x}f(g)d\lambda^x(g)$. We denote by $\int f d\lambda$ the element of $C_0(G^{(0)})$ obtained when integrating a compactly supported function $f$ w.r.t. the Haar system. If $A$ is a $G$-algebra and $E$ is a $A$-Hilbert module, $f\mapsto \int f d\lambda$ extends by linearity and continuity to $C_c(G,A)\subseteq r^*A $ and $C_c(G,E)\subseteq r^*E$. We will still denote by $\int fd\lambda$ the element in $A$ or $E$ if $f$ is in $r^* A$ or $r^* E$ respectively. %Notice that the invariance of the Haar system implies that, if $V\in\mathcal L_{s^*B}(s^*E,r^*E)$ denotes the action of $G$ on $E$, then $\int V F V^* d\lambda = \int F d\lambda$ in $\mathcal L_B(E)$ , for every $F\in\mathcal L_B(E)$. 

\begin{lem}[lemma $3.6$,\cite{TuBC2}]\label{JLTform}
Let $X$ be a $G$-compact proper $G$-space such that the anchor map $p:X\rightarrow G^{(0)}$ is locally injective, and let $B$ be a $G$-algebra. Then for every $z\in RK^G(X,B)$ there exists a $G$-proper $G$-compact space $Z\subseteq X$ and a $K$-cycle $(H_B, \pi, T)\in \mathbb E^G(C_0(Z),B)$ representing $z$ such that :
\begin{itemize}
\item[$\bullet$] $T$ is self-adjoint and $-1 \leq T\leq 1$,
\item[$\bullet$] $T$ is $G$-equivariant, i.e. $r^* T = V s^*T V^*$ ,
\item[$\bullet$] $T$ commutes with the action of $Z$, i.e. $[\pi(a),T]= 0$ for all $a\in C_0(Z)$.
\end{itemize}
\end{lem}

%%%%%%%%%%%%% NOUVELLE PREUVE
\begin{dem}
Let $(E,\pi,T)\in \mathbb E(C_0(X),B)$ be a $K$-cycle. Denote by $\alpha : s^* C_0(Z) \rightarrow r^*C_0(Z)$ the action of $G$. Let $K$ be a compact fundamental domain for the action of $G$ on $Z$. By local injectivity of $p$, let $(U_j)_j$ be a finite open cover of $K$ such that $p_{|U_j}$ is injective for every $U_j\in \mathcal U$. There exist compactly supported continuous functions $\phi_j : Z\rightarrow [0,1]$ such that  \[\text{supp }\phi_j \subseteq U_j \quad\text{ and }\quad K\subseteq \cup_{j} \phi_j^{-1}(0,+\infty).\]
Up to replacing $\phi_j$ by the continuous function $\phi_j(z) / \sum_{k,g} \phi_k(z.g)$, we can assume 
\[\sum_{j,g\in G^{p(z)}} \phi_U (z.g) = 1,\forall z\in Z.\] 
The latter is indeed continuous, as it can be expressed as $\frac{\phi_j} {(\sum_j\int \phi_j d\lambda ) }$. The condition $\sum_{j,g\in G^{p(z)}} \phi_j (z.g) = 1,\forall z\in Z$ implies that 
\[\sum_j \int V(s^*\pi(\phi_j))V^*d\lambda = id_E ,\]
in the sense of a weak integral, i.e. the equality holds when evaluated on elements of $E$.\\
  
Indeed, $\pi$ is $G$-equivariant, and by invariance of the Haar system, $\sum_j\int \alpha(s^*\phi_j)d\lambda = 1_{C_0(Z)}$. Composing with $\pi$ gives the equality by continuity and linearity of $\pi$.\\  

Define $F_j = \pi(\phi_j^\frac{1}{2}) F \pi(\phi_j^\frac{1}{2})\in\mathcal L_B(E)$. Then $s^* F_j\in C_c(G,\mathcal L_{s^* B}(s^* E))$ and set
\[F'= \sum_{j} \int V (s^*F_j) V^* d\lambda.\] 
By invariance of the Haar system, the operator $F'$ is  $G$-invariant. Moreover it commutes with the action of $C_0(Z)$. Indeed, we can compute its fibers : for every $x\in G^{(0)}$,
\[F'_x = \sum_{j, g\in G^x} V_g\pi(\phi_j^{\frac{1}{2}}) F_{s(g)}\pi(\phi_j^{\frac{1}{2}})V_g^*.\] 
Hence, by local injectivity, for all $g\in G$, there exists $z_g\in Z$ such that $Z_{s(g)}\cap U = \{z_g\}$, and for all $f\in C_0(Z_{s(g)})$, we have that $\phi_U^{\frac{1}{2}} f = f(z_g) \phi_U^{\frac{1}{2}}$, hence $[ \phi_U^{\frac{1}{2}} T_{s(g)} \phi_U^{\frac{1}{2}},f ] = 0 $.\\

Moreover, $F'$ is a compact perturbation of $F$ as the following computation shows.\\
\[\begin{array}{rl}
F -F' 	& = (\sum_{j} \int V (s^* \phi_j) V^*d\lambda)F -F' \\
		& = \sum_{j}\int V (s^*\phi_j^{\frac{1}{2}}) \ 
\left( (s^*\phi_j^\frac{1}{2}) V^* (r^*F) V- (s^*F)(s^*\phi_j^\frac{1}{2}) \right) \ V^* \ d\lambda\\
		& = \sum_{j} \int V (s^*\phi_j^{\frac{1}{2}}) \ 
			\left( (s^*\phi_j^\frac{1}{2}) (V^* (r^*F)V - (s^*F)) + s^*(\phi_j^\frac{1}{2} [\phi_j^\frac{1}{2},F] )\right) \ 
				V^* \ d\lambda	  	
\end{array}\]
Each of the summand is compact, hence $[H_B, \pi,T]=[H_B,\pi,T']$.\\ 
\qed
\end{dem}
%%%%%%%%%%%%%%%%%%%%

\begin{lem}
Let $X$ be a $G$-compact proper $G$-space such that the anchor map $p:X\rightarrow G^{(0)}$ is locally injective, and let $(B_j)_j$ be a countable family of $G$-algebras. Then the projection $\prod_j B_j \otimes \mathfrak K \rightarrow B_j\otimes \mathfrak K$ induces an isomorphism
\[\Theta : RK^G(X,\prod_j B_j\otimes \mathfrak K)\rightarrow \prod_j RK^G(X,B_j\otimes \mathfrak K)\cong \prod_j RK^G(X,B_j).\]
\label{LocalInjectivity}
\end{lem}

% NEW PROOF
\begin{dem}
Set $B_\infty = \prod_j B_j\otimes \mathfrak K $. Let us define a $Z_2$-graded homomorphism 
\[\eta :  \prod_j RK^G(X,B_j\otimes \mathfrak K) \rightarrow RK^G(X,\prod_j B_j\otimes \mathfrak K).\] 
Let $Z\subseteq X$ be a $G$-proper $G$-compact subspace, and, for all $j$, let $z_j\in KK^G(C_0(Z),B_j)$ be represented by a standard $K$-cycle $(H_{B_j},\pi_j ,T_j)$. By lemma \ref{JLTform}, we can suppose that $T_j$ is $G$-equivariant self-adjoint, commutes with the action of $C_0(Z)$ and $-1\leq T_j \leq 1$. \\

Define $E_\infty = \prod_j (E_j\otimes \mathfrak K)$. It is a $B_\infty$-Hilbert module with respect to the scalar product 
\[\langle \xi,\eta\rangle = \prod (\langle \xi_j,\eta_j\rangle \otimes F_j^* F'_j )\in B_\infty\] 
for every $\xi = (\xi_j\otimes F_j)$ and $\eta = (\eta_j\otimes F'_j)$ in $E_\infty$. Define $T = (T_j\otimes E_{11} + id_{E_j}\otimes (1-E_{11}))_j$, where $E_{11}(x)= \langle e_1,x\rangle. e_1$ is the rank-one operator projecting on $e_1$. ($\{e_j\}$ denote the canonical orthonormal basis of $l^2(\mathbb N)$).
These conditions ensure that $T $ defines an operator in $\mathcal L_{B_\infty}(E_\infty)$.\\

Define $\pi(a) = (\pi_j(a)\otimes id_{\mathfrak K})_j\in  \mathcal L_{B_\infty}(E_\infty)$, as $\sup ||\pi_j(a)||\leq ||a||$. Then, $[\pi(a),T]=0$, $T_j^*=T_j$ and $r^*T = Vs^*T V^*$. %Moreover $0\leq T_j^2-1\leq 1$, hence, according to the remark \ref{construction},
Let us show that $T^2-1$ is in $\mathfrak K_{B_\infty}(E_\infty)$.\\ 

Indeed, each $F_j = T_j^2-1\in \mathfrak K \otimes B_j$ can be approximated by a finite rank operator, hence there exist $N_j>0$, and $\xi_k^{(j)},\eta_k^{(j)}\in H_{B_j}$ such that $F_j= \sum_{k=1}^{N_j}\theta_{\xi_k^{(j)},\eta_k^{(j)}}$ for every $j$. A computation shows that $(F_j\otimes E_{11})_j = \theta_{\sum_k \tilde\xi_k,\sum_k \tilde \eta_k}$ where $\tilde \xi_k = (\xi_k^{(j)}\otimes E_{j1})$ and $\tilde \eta_k = (\eta_k^{(j)}\otimes E_{j1})$, which is a rank one operator, hence $T^2-1$ is compact since a product of rank one operators is of rank one.\\ 
%Indeed, these elements are uniformly bounded elements of $\prod \mathfrak K_{B_j\otimes\mathfrak K}$, which are in $\mathfrak K_{B_\infty}(H_{B_\infty})$. 
This ensures that $[E_\infty,\pi,T]\in KK^G(C_0(Z),B_\infty)$. Let $\eta((z_j)_j) = [H_{B_\infty},\pi,T] $. It is clear that $\Theta\circ \eta = id_{\prod_j RK^G(X,B_j\otimes K) }$, and $\Theta $ is onto.\\

Let us prove that $\Theta $ is one to one. Let $z\in KK^G(C_0(Z),B_\infty)$ such that $\Theta(z) = 0$. Let us denote $\Theta(z) = (z_j)_j$ where each $z_j\in KK^G(C_0(Z),B_j)$ is represented by $K$-cycles $(H_{B_j}, \pi_j,T_j)$ homotopic to a degenerate $K$-cycle by an operator homotopy 
\[\left\{\begin{array}{rcl} [0,1] & \rightarrow & \mathbb E^G(C_0(Z),B_\infty)\\ s & \mapsto & (H_{B_j}, \pi_j,T_j(s))\end{array}\right.\] 
such that $T_j(s)$ is $G$-equivariant self-adjoint, commutes with the action of $C_0(Z)$ and $-1\leq T_j \leq 1$ for every $s\in [0,1]$.\\

%%%
Set 
\[\mathcal{\tilde C}_j = \{T \in\mathcal L_{B_j}(H_{B_j}) \text{ s.t. } [\pi(a),T]=0 \ \forall a\in C_0(Z)\},\] 
and let $\mathcal C_j$ be the closed ideal $\{ T \in \mathcal {\tilde C}_j \text{ s.t. } \pi(a)T\in {\mathfrak K}_{B_j}(H_{B_j})  \ \forall a\in C_0(Z)\}$. Similarly for $\mathcal L_{B_\infty}(E_\infty)$ and $\mathfrak K_{B_\infty}(E_\infty)$, define $\mathcal{\tilde C}$ and the closed ideal $\mathcal{C}$. Our goal is to show that the family of operator homotopies can be lifted to a global one.\\

Denote by $\overline T_j$ the class of $T_j$ in $\mathcal{\tilde C}_j / \mathcal C_j $. For every $j$, let $t_0=0< t_1 < ... < t_{l_j}=1$ be a partition of $[0,1]$ such that $||\overline{T}_j(t_{k+1})-\overline{T}_j(t_k)||< 1$ for every $k\in\{0,..,l_j-1\}$. Put $T_j(t_k) = T_{j,k}$. Then 
%Denote by $\overline T_j$ the class of $T_j$ in $\mathcal L_{B_j}(H_{B_j})/ {\mathfrak K}_{B_j}(H_{B_j}) $. For every $j$, let $t_0=0< t_1 < ... < t_l=1$ be a partition of $[0,1]$ such that $||\overline{T}_j(t_{k+1})-\overline{T}_j(t_k)||< 1$ for every $k\in\{0,..,l-1\}$. Put $T_j(t_k) = T_{j,k}$. Then 
\[diag(\overline{T}_{j,0},1,1,1,1,...)= diag(1,\overline{T}_{j,1}^*,...,\overline{T}_{j,l}^*,1,1,1,...) \ diag(\overline{T}_{j,0},\overline{T}_{j,1},...,\overline{T}_{j,l},1,1,1,...)\]
is $2$-Lipschitz homotopic to 
\[diag(\overline{T}_{j,1}^*,...,\overline{T}_{j,l}^*,1,1,1,1,...) \ diag(\overline{T}_{j,0},\overline{T}_{j,1},...,\overline{T}_{j,l},1,1,1,...)\]\[=diag(\overline{T}_{j,1}^*\overline{T}_{j,0} ,..., \overline{T}_{j,l}^*\overline{T}_{j,l-1},\overline{T}_{j,l},1,1,1,...) \] %diag(e^{ a_j^{(1)}} ,..., e^{ a_j^{(l)}},1,1,1,...)
by rotations. We obtained that \[diag(\overline{T}_{j,0},1,1,1,1,...)\] 
is $2$-Lipschitz homotopic to $diag(\overline{T}_{j,l},\overline{T}_{j,1}^*\overline{T}_{j,0} ,..., \overline{T}_{j,l}^*\overline{T}_{j,l-1},1,1,1,...)$.\\

By (\cite{WeggeOlsen}, Proposition $4.2.4$) for every $j$ and $k$, there exists $a_j^{(k)}\in \mathcal{\tilde C}_j / \mathcal C_j$ such that : 
%$a_j^{(k)}\in \mathcal L_{B_j}(H_{B_j})/ {\mathfrak K}_{B_j}(H_{B_j})$ such that :
\[\overline{T}_{j,k+1}^*\overline{T}_{j,k}=e^{ a_j^{(k)}} \text{ and such that } ||a_j^{(k)}||<1.\] 
%Then (\cite{WeggeOlsen}, Lemma $17.3.3$) we can find a lift $x_j^{(k)}\in  \mathcal L_{B_j}(H_{B_j})$ of $a_j^{(k)}$ such that $||a_j^{(k)}||= ||x_j^{(k)}||$ for every $j$ and $k$.
%For each $k$, choose a continuous function $\phi_k : [0,1]\rightarrow [0,1]$ such that :\\
%\begin{itemize}
%\item[$\bullet$] $\phi_j $ is constant equal to $0$ on $[0,t_k ]$,
%\item[$\bullet$] $\phi_j$ is affine on $[t_k,t_{k+1}]$,
%\item[$\bullet$] $\phi_j $ is constant equal to $1$ on $[t_{k+1},1 ]$.\\
%\end{itemize}
Define 
\[\tilde T_j(s) = \begin{pmatrix} \overline{T}_{j,l} & 0 \\ 0 & \exp ( s \ diag(a_j^{(1)},..., a_j^{(l)},0,0,0,...))\end{pmatrix}.\]
%\[\tilde T_j(s) = diag(e^{\phi_j(s) a_j^{(1)}},..., e^{\phi_j(s) a_j^{(l)}},\overline{T}_{j,l},1,1,1,...).\] 
Then $s\mapsto \tilde T_j(s)$ composed with the first homotopy by rotation and permutation is $L$-Lipschitz for every $j$ for some constant $L$ independent of $j$, and it provides an %homotopy between $diag(\overline{T}_j,1,1,..)$ and a degenerate operator. Hence we get a $L$-Lipschitz 
element of $C([0,1],\mathcal{\tilde C}/\mathcal{C})$. We can lift it to $T\in C([0,1],\mathcal{\tilde C})$, which gives a homotopy 
\[\left\{\begin{array}{rcl} [0,1] & \rightarrow & \mathbb E^G(C_0(Z),B_\infty)\\ s & \mapsto & (E_\infty, \pi,T(s))\end{array}\right.\]
%defines a $L$-Lipschitz operator homotopy 
between $0$ and $z$.\\

\qed
\end{dem}
%END NEW PROOF

\begin{lem}\label{prod}
Let $G$ be a locally compact, $\sigma$-compact étale groupoid, $\{B_j\}_{j\geq  0}$ a family of $G$-algebras and $\mathfrak K$ the algebra of compact operators over a separable Hilbert space. Then, for every finite dimensional proper $G$-compact $G$-simplicial complex $\Delta$, we have an $\Z_2$-graded isomorphism of abelian groups
\[RK^G(\Delta,\prod_j B_j\otimes \mathfrak K)\cong \prod_j RK^G(\Delta,B_j)\]
\end{lem}

\begin{dem}
For all $j$ and any locally compact $G$-space $X$, the projection $\prod_j B_j\otimes \mathfrak K\rightarrow B_j \otimes \mathfrak K$ induces a morphism
\[\Theta^X : KK^G(C_0(X),\prod_j B_j\otimes \mathfrak K )\rightarrow \prod_j  KK^G(C_0(X),B_j\otimes \mathfrak K ).\]
Let $X_0\subseteq X_1 \subseteq...\subseteq X_n$ be the $n$-skeleton decomposition associated to the simplicial structure of typed $G$-simplicial complex $\Delta$ and let $Z_j = C_0(X_j)$, $Z^j_{j-1} = C_0(X_j \setminus X_{j-1})$ and $\Theta_j = \Theta^{X_j}$.
We will show the claim by induction on the dimension of $\Delta$.\\

By naturality of the boundary element, the extension of $G$-algebras $0\rightarrow Z^j_{j-1} \rightarrow Z_j \rightarrow Z_{j-1}\rightarrow 0$ gives a commutative diagram with exact rows :
%lines :
%\[\begin{tikzcd}
%KK_*(Z^j_{j-1},\prod_j B_j\otimes K)\arrow{r}{\delta}\arrow{d}{\Theta^j_{j-1}} & KK_*(Z_{j-1},\prod_j B_j\otimes K)\arrow{r}\arrow{d}{\Theta_{j-1}} & KK_*(Z_j,\prod_j B_j\otimes K)\arrow{r}\arrow{d}{\Theta_j} & KK_*(Z^j_{j-1},\prod_j B_j\otimes K)\arrow{r}{\delta} \arrow{d}{\Theta^j_{j-1}} & KK_*(Z_{j-1},\prod_j B_j\otimes K)\arrow{d}{\Theta_{j-1}}\\
%\prod_j KK_*(\tilde Z^j_{j-1},B_j \otimes K)\arrow{r}{\delta} & \prod_j KK_*(\tilde Z_{j-1},B_j \otimes K)\arrow{r} & \prod_j KK_*(\tilde Z_j,B_j \otimes K)\arrow{r} & \prod_j KK_*(\tilde Z^j_{j-1},B_j \otimes K)\arrow{r}{\delta} & \prod_j KK_*(\tilde Z_{j-1},B_j \otimes K)\\
%\end{tikzcd}\]
\[\begin{tikzcd}
KK^G_*(Z^j_{j-1},\prod_j B_j\otimes \mathfrak K)\arrow{d}{\partial}\arrow{r}{\Theta^j_{j-1}} & \prod_j KK^G_*( Z^j_{j-1},B_j \otimes \mathfrak K)\arrow{d}{\partial}  \\
KK^G_*(Z_{j-1},\prod_j B_j\otimes \mathfrak K)\arrow{d}\arrow{r}{\Theta_{j-1}}  & \prod_j KK^G_*( Z_{j-1},B_j \otimes \mathfrak K)\arrow{d} \\
KK^G_*(Z_j,\prod_j B_j\otimes \mathfrak K)\arrow{d}\arrow{r}{\Theta_j} & \prod_j KK^G_*( Z_j,B_j \otimes \mathfrak K)\arrow{d} \\
KK^G_*(Z^j_{j-1},\prod_j B_j\otimes \mathfrak K)\arrow{d}{\partial}\arrow{r}{\Theta^j_{j-1}} & \prod_j KK^G_*( Z^j_{j-1},B_j \otimes \mathfrak K)\arrow{d}{\partial}\\
KK^G_*(Z_{j-1},\prod_j B_j\otimes \mathfrak K)\arrow{r}{\Theta_{j-1}} & \prod_j KK^G_*( Z_{j-1},B_j \otimes \mathfrak K)
\end{tikzcd}\]

The five lemma ensures that if $\Theta_{j-1}$ and $\Theta^j_{j-1}$ are isomorphisms, then so is $\Theta_j$. Moreover, because $\Delta$ is a typed $G$-simplicial complex (see \ref{Gcomplex}), $X_j\setminus X_{j-1}$ is $G$-equivariantly homeomorphic to $\mathring \sigma_j \times \Sigma_j$, where  $\mathring \sigma _j $ denotes the interior of the standard simplex, $\Sigma_j$ is the set of centers of $j$-simplices of $X_j$, and where $G$ acts trivially on $\mathring \sigma _j$. Bott periodicity ensures then that the following diagram commutes :
\[\begin{tikzcd}
KK^G(Z_{j-1}^j,\prod_k B_k\otimes \mathfrak K) \arrow{r}{\Theta_{j-1}^j}\arrow{d}{\cong} & \prod_k KK^G(Z_{j-1}^j,B_k) \arrow{d}{\cong}\\
KK^G(\Sigma_j,\prod_k B_k\otimes \mathfrak K) \arrow{r}{\Theta^{\Sigma_j}} & \prod_k KK^G(\Sigma_j ,B_k)\\ 
\end{tikzcd}\]
with vertical arrows being isomorphisms given by Bott periodicity. By lemma \ref{LocalInjectivity}, $\Theta^{\Sigma_j}$ is an isomorphism, hence $\Theta_{j-1}^j$ is an isomorphism. We proved that if $\Theta_{j-1}$ is an isomorphism, then so is $\Theta_j$. By induction, proving that $\Theta_0$ is an isomorphism concludes the proof, which is essentially the content of lemma \ref{LocalInjectivity} : $X_0$ is a $G$-compact proper $G$-space, and its anchor map is just the target map $r:G\rightarrow G^{(0)}$, which is supposed to be étale, so locally injective.\\
\qed
\end{dem}

%%%%%%%% CONTROLLED K THEORY AND PRODUCTS

We now turn our attention to the behaviour of controlled $K$-theory with respect to products. Let $\mathcal B = (B_j)_j$ a countable family of $\mathcal E$-filtered $C^*$-algebras. For every $E\in\mathcal E$, put 
\[\mathcal B_E = \prod_{j} (B_j)_E\otimes \mathfrak K.\]
Define $\mathcal B_\infty$ as the closure of $\cup_{E\in\mathcal E}\mathcal B_E$ in $\prod_j (B_j\otimes \mathfrak K)$. This construction is called the controlled product of a family of filtered $C^*$-algebras, and is also a $\mathcal E$-filtered $C^*$-algebra. The following lemma (\cite{OY3}, Lemma $1.14$) relates the controlled $K$-theory of $\mathcal B_\infty$ and the product of the controlled $K$-theory groups of $B_j$.\\

\begin{lem}[\cite{OY3}] Let $\mathcal B$ be a countable family of $\mathcal E$-filtered $C^*$-algebras. With the previous notations, and for every $\varepsilon\in(0,\frac{1}{4})$ and every $E\in\mathcal E$, there exists a control pair $(\alpha, h)$, independent of the family $\mathcal B$, such that the family of  maps 
\[K^{\varepsilon,E}(\mathcal B_\infty)\rightarrow \prod K^{\varepsilon,E}(B_j)\]
induced by the composition $\prod B_j\otimes \mathfrak K\rightarrow B_k\otimes \mathfrak K$ and the controlled Morita equivalence gives a $(\alpha,h)$-controlled isomorphism
\[\hat K(\mathcal B_\infty)\rightarrow \prod \hat K(B_j).\]
\end{lem}

\begin{lem} Let $G$ be an étale groupoid with compact base space $G^{(0)}$. Let $(B_j)_j$ be a countable family of $G$-algebras, and $\mathfrak K$ the $G$-algebra of compact operators on a separable Hilbert space with trivial $G$-action. Then, $\mathcal B^G = ((B_j\otimes \mathfrak K)\rtimes_r G)_j$ is a countable family of $\mathcal E$-filtered $C^*$-algebra, and there exists a filtered $*$-isomorphism :
\[ (\prod_j B_j\otimes\mathfrak K)\rtimes_r G \rightarrow  \mathcal B^G_\infty.\]
\end{lem}

\begin{dem} Notice that, as $\mathfrak K$ is endowed with the trivial action, $(B_j\rtimes_r G)\otimes\mathfrak K \cong (B_j\otimes\mathfrak K)\rtimes_r G$ for every $j$. By definition, for every $E\in \mathcal E$,  $\left( (\prod B_j\otimes\mathfrak K)\rtimes_r G \right)_E = \prod_j \left( (B_j\otimes\mathfrak K)\rtimes_r G\right)_E \cong \prod_j \left((B_j\rtimes_r G)_E \otimes \mathfrak K \right) = (\mathcal B^G_\infty )_E$.\\
\qed
\end{dem}

The two previous lemmas entail the following result.

\begin{cor} \label{controlledprod} Let $G$ be an étale groupoid with compact base space $G^{(0)}$. Let $(B_j)_j$ be a countable family of $G$-algebras, and $\mathfrak K$ the $G$-algebra of compact operators on a separable Hilbert space with trivial $G$-action. Then, there exists a control pair $(\alpha, h)$ and a $(\alpha, h)$-controlled isomorphism
\[ \hat K(\mathcal B^G_\infty) \rightarrow \prod \hat K(B_j\rtimes_r G). \]
\end{cor}

%%%%%%%%%%%%%%%%%%%%%%%%%%%%%%%%%%%%%%
\subsection{Quantitative statements}
%%%%%%%%%%%%%%%%%%%%%%%%%%%%%%%%%%%%%%

We can now prove the quantitative statements.

\begin{thm}\label{Quant1}
Let $B$ a $G$-algebra, and $\tilde B = l^\infty(\N,B\otimes \mathfrak K)$. Then $\mu_{G,\tilde B}$ is injective if and only if for every $E\in\mathcal E,\varepsilon\in(0,\frac{1}{4})$ and $F$ such that $k_J(\varepsilon).E\subseteq F$, there exists $E' \in\mathcal E$ such that $E\subseteq E'$ and $QI_{G,B}(E,E',\varepsilon,F)$ holds. 
\end{thm}

\begin{dem}
Let $x\in RK^G(P_E(G),\tilde B)$ such that $\mu_{G,\tilde B}^E(x)=0$. Then, as the quantitative assembly maps factorize $\mu_{G,\tilde B}$, there exist $\varepsilon>0$ and $F$ such that $k_J(\varepsilon). E \subseteq F$, satisfying $\mu_{G,\tilde B}^{\varepsilon,E,F}(x)=0$. Let us denote by $(x_j)_{j\in \N}$ the element of $\prod_j RK^G(P_E(G),A)$ corresponding to $x$ under the isomorphism of lemma \ref{prod}. Now let $F'$ be in $\mathcal E$ such that $F\subseteq F'$ and such that $QI_{A}(E,E',\varepsilon,F)$ holds. That ensures that $q_E^{E'}(x_j)=0$ in $RK^G(P_{E'}(G),B)$, and $q_E^{E'}(x)=0$ hence $\mu_{G,\tilde B}$ is injective.\\

% Using the isomorphism of lemma \ref{prod} and the Morita equivalence, we can identify $x$ with $(x_j)_j$ under $RK^G(P_E(G),\tilde A)\simeq\prod_j RK^G(P_E(G),A)$. Now let $F'$ such that $F\subseteq F'$ and $QI_{A}(E,E',\varepsilon,F)$ holds. That ensures that $x_j=0$ in $RK^G(P_{E'}(G),B)$, and $x=0$.\\

For the converse, suppose one can find $E,\varepsilon,F$ such that $QI_{G,B}(E,E',\varepsilon,F)$ is NOT true for all $F'$ such that $F\subseteq F'$. Then, by $\sigma$-compactness, one can extract a increasing exhausting sequence $E_j\in \mathcal E$ such that $\cup E_j =G$, $E\subseteq E_j$ and $x_j\in RK^G(P_E(G),B)$ such that $\mu_{G,\tilde B}^{\varepsilon,E,F}(x_j)=0$ and $q_E^{E_j}(x_j)\neq 0$ in $RK^G(P_{E_j}(G),B)$. Let $x$ be the image of $(x_j)\in \prod RK^G(P_E(G),B)$ in $RK^G(P_E(G),\tilde B)$ under the isomorphism of lemma \ref{prod}. By corollary \ref{controlledprod}, and up to a rescaling, we have $\mu_{G,\tilde B}^E(x)=0$, and $q_E^{E'}(x)\neq 0$ in $RK^G(P_{E'}(G),\tilde A)$ for at least one $E'$ such that $E\subseteq E'$ by exhaustivity, so $\mu_{G,\tilde B}$ is not injective. \\
\qed   
\end{dem}

We also have a theorem relating quantitative surjectivity for $\hat\mu_{G,B}$ and surjectivity of $\mu_{G,\tilde B}$.

\begin{thm}\label{Quant2}
Let $B$ a $G$-algebra, and $\tilde B = l^\infty(\N,B\otimes \mathfrak K)$. Then there exists $\lambda \geq 1$ such that $\mu_{G,\tilde B}$ is onto if and only if for any $\varepsilon\in ( 0 ,\frac{1}{4\lambda})$ and every $F\in\mathcal E$, there exist $E,F'\in\mathcal E$ such that  $k_J(\varepsilon) .E \subseteq F$, $F\subseteq F'$ and such that $QS_{B,G}(E,F,F',\varepsilon,\lambda\varepsilon)$ holds.
\end{thm}

\begin{dem}
Let $\lambda\geq 1$ the universal constant of remark \ref{approximation} : for any $C^*$-algebra and any $x\in K^{\varepsilon, F}(A)$ such that $\iota_{\varepsilon,F} (x) =0$, there exists $F'$ such that $F\subseteq F'$ and $\iota_{\varepsilon,F}^{\lambda\varepsilon,F'} (x) =0$.\\

Let $y\in K_*(\tilde B\rtimes_r G)$. By remark \ref{approximation}, there exist $\varepsilon\in (0,\frac{1}{4})$, $F\in\mathcal E$ and $z\in K^{\varepsilon,F}(\tilde B\rtimes G)$ such that $\iota_{\varepsilon,F}(z) = y$. Up to a rescaling of the parameters, let $(z_j)$ the element of $\prod_j K^{\varepsilon,F}(B\rtimes_r G)$ corresponding to $z$ under the controlled isomorphism of corollary \ref{controlledprod}. Let $E$ and $F'$ such that $k_J(\varepsilon).E\subseteq F$ and $QS(E,F,F',\varepsilon,\lambda\varepsilon)$ : for every $j$, there exists $x_j\in RK^G(P_E(G),B) $ such that $\mu_{G,B}^{\lambda\varepsilon,E,F'}(x_j)=\iota_{\varepsilon,F}^{\lambda\varepsilon,F'}(z_j)$. 
Let $x\in RK^G(P_E(G),\tilde B)$ be the element of $RK^G(P_E(G),\tilde B)$ corresponding to $(x_j)\in \prod_j RK^G(P_E(G),B)$ under the isomorphism of lemma \ref{prod}. Naturality of the assembly maps, and compatibility of the controlled assembly maps with the usual one ensures that $\mu_{G,\tilde B}^E(x)=z$, whereby $\mu_{G,\tilde B}$ is onto.\\

Assume that there exist $\varepsilon \in (0,\frac{1}{4\lambda})$ and a nonempty controlled subset $F\in\mathcal E$ such that for every $E,F'\in\mathcal E$ such that $k_J(\varepsilon). E \subseteq F$, $QS_{G,B}(E,F,F',\varepsilon,\lambda\varepsilon)$ does not hold. Let $(E_j)$ and $(F_j)$ be unbounded increasing sequences of controlled subsets such that $(E_j)$ is an exhausting family, $F\subseteq F_j$ and $k_J(\varepsilon).E_j\subseteq F_j$. Let $y_j\in K^{\varepsilon,F}(B\rtimes_r G)$ such that $\iota_{\varepsilon,F}^{\lambda\varepsilon,F_j}(y_j)$ is not in the range of $\mu_{G,B}^{\lambda\varepsilon, E_j , F_j}$. Let $y\in K^{\varepsilon,F}(\tilde B\rtimes_r G)$ be the element of corresponding to $(y_j)_j$ under the isomorphism of corollary \ref{controlledprod}, up to a rescaling of the parameters. If there exists $x\in RK^G(P_{E'}(G),\tilde B)$ for a $E'\in\mathcal E$ such that $E\subseteq E'$ and $\iota_{\varepsilon,F}( y) =\mu_{G,\tilde B}^{E'}(x)$ then there would exists a $F'\in\mathcal E$ such that $F\subseteq F'$ and
\[\iota_{\varepsilon,F}^{\lambda\varepsilon, F'}( y) =\mu_{G,\tilde B}^{\lambda\varepsilon,E',F'}(x) = \iota_{\varepsilon,F}^{\lambda\varepsilon, F'}\circ\mu_{G,\tilde B}^{\varepsilon,E',F}(x).\]
Now choose $j$ such that $E'\subseteq E_j $ and $F'\subseteq F_j$, and compose the previous equality with $\iota_{\lambda\varepsilon,F'}^{\lambda\varepsilon,F_j}$ and $q_{E'}^{E_j}$ to obtain $\iota_{\varepsilon,F}^{\lambda\varepsilon,F_j}(y_j)=\mu_{G,B}^{\lambda\varepsilon , E_j,F_j}(x_j)$ which contradicts our assumption. Hence $\mu_{G,\tilde B}$ is not onto. \\
\qed
\end{dem}

Using the only if part of the proofs of the quantitative statements, we can easily prove the following theorem by replacing $l^\infty(\N,B\otimes \mathfrak K)$ by $\prod (B_j \otimes \mathfrak K)$.

\begin{thm}\label{UniformQS} Let $G$ be an étale groupoid with compact base space. 
\begin{itemize}
\item[$\bullet$] Assume that for any $G$-algebra $B$, $\mu_{G,B}$ is one-to-one. Then, for any $\varepsilon\in (0,\frac{1}{4})$ and every $E,F\in\mathcal E$ such that $k_J(\varepsilon). E\subseteq F$, there exists $E'\in\mathcal E$ such that $E\subseteq E'$ and such that $QI_{G,A}(E,E',\varepsilon,F)$ holds for any $G$-algebra $B$.
\item[$\bullet$] Assume that for any $G$-algebra $B$, $\mu_{G,B}$ is onto. Then, for some $\lambda\geq 1$ and for any $\varepsilon\in (0,\frac{1}{4\lambda})$ and every $F\in\mathcal E$, there exists $E,F'\in\mathcal E$ such that $k_J(\varepsilon). E\subseteq F'$ and $F\subseteq F'$ such that, for any $G$-algebra $B$, $QS_{G,A}(E, F,F',\varepsilon,\lambda \varepsilon)$ holds.
\end{itemize}
\end{thm}

The Theorems \ref{Quant1} and \ref{Quant2} imply together Theorem \ref{ControlledBC}.

%%%%%%%%%%%%%%%%%%%%%%%%%%%%%%%%%%%%%%%%%%%%%%%
\subsection{Persistence approximation property}
%%%%%%%%%%%%%%%%%%%%%%%%%%%%%%%%%%%%%%%%%%%%%%%

We recall the following definition from \cite{OY3}.

\begin{definition}
Let $B$ be a $\mathcal E$-filtered $C^*$-algebra, $\lambda>0$, $\varepsilon,\varepsilon'$ be positive numbers such that $0<\varepsilon <\varepsilon' <\frac{1}{4}$ and $F,F'\in\mathcal E$ be nonempty controlled subsets such that $F\subseteq F'$. The following property is called Persistance Approximation Property :
\begin{itemize}
\item[$\bullet$] $PA_B(\varepsilon,\varepsilon',F,F')$ : for every $x\in K_*^{\varepsilon,F}(B)$ such that $\iota_{\varepsilon,F}(x)=0$ in $K_*(B)$, then $\iota_{\varepsilon,F}^{\varepsilon',F'}(x)=0$ in $K_*^{\varepsilon',F'}(B)$.
\item[$\bullet$] $B$ is said to satisfy the Persistance Approximation Property $(PAP)_\lambda$ if for every nonempty $F\in\mathcal E$ and every $\varepsilon\in (0,\frac{1}{4\lambda})$, there exists $F'\in\mathcal E$ nonempty such that $PA_B(\varepsilon,\lambda\varepsilon,F,F')$ holds.
\end{itemize}
\end{definition}

%\begin{thm} \label{PAPG}
%Up to some hypothesis, if $\mu_{G,l^\infty(\N, K_A)}$ is onto and $\mu_{G,A}$ is one-to-one, then, for a universal constant $\lambda_{PA}$, for all $\varepsilon \in(0,\frac{1}{4\lambda_{PA}})$ and nonempty $F$, there exists $F'$ such that $F\subseteq F'$ and $PA_{A\rtimes G}(\varepsilon,\lambda_{PA}\varepsilon,F,F')$ holds.
%\end{thm}

The following result gives a sufficient condition for $(PAP)_\lambda$ to be satisfied for a large class of $C^*$-algebras.

\begin{thm} \label{PAPG}
Let $G$ be an étale groupoid such that :
\begin{itemize}
\item[$\bullet$] $G^{(0)}$ is compact, 
\item[$\bullet$] $G$ admits a cocompact example for universal space for proper actions. 
\end{itemize}
Then there exists a universal constant $\lambda_{PA}\geq 1$ such that, for every $G$-algebra $A$, if $\mu_{G,l^\infty(\N, A\otimes\mathfrak K)}$ is onto and $\mu_{G,A}$ is one-to-one, then for every $\varepsilon \in(0,\frac{1}{4\lambda_{PA}})$ and every nonempty $F\in\mathcal E$, there exists $F'\in\mathcal E$ such that $F\subseteq F'$ and $PA_{A\rtimes_r G}(\varepsilon,\lambda_{PA}\varepsilon,F,F')$ holds.\\
\end{thm}

\begin{dem}
Let us denote $l^\infty(\N,A\otimes\mathfrak K)$ by $\tilde A$. Let $(\alpha,h)$ be the control pair of the controlled isomorphism of corollary \ref{controlledprod}, and $\lambda\geq 1$ be the constant of remark \ref{approximation}. Set $\lambda_{PA} = \lambda\alpha$. \\

Assume the statement does not holds : there exists $\varepsilon$ and $F$ such that $PA_{A\rtimes_r G}(\varepsilon,\lambda_{PA}\varepsilon,F,F')$ is not true for every $F'\in\mathcal E$ such that $F \subseteq F'$. Then we can extract an increasing exhausting sequence of controlled subsets $F_j$ such that $F\subseteq F_j$ and elements $x_j\in K^{\varepsilon,F}(A\rtimes_r G)$ such that $\iota_{\varepsilon,F}(x_j)=0$ and $\iota_{\varepsilon,F}^{\lambda_{PA}\varepsilon,F_j}(x_j)\neq 0$. \\

Let $x$ be the element of $K^{\alpha\varepsilon,h_\varepsilon  F}(\tilde A\rtimes_r G)$ corresponding to \[(x_j)\in\prod_j K^{\varepsilon,F}(A\rtimes_r G)\] under the controlled isomorphism of corollary \ref{controlledprod}. Recall that the following diagram commutes
%According to the lemma \ref{LocalInjectivity}, there exists $x\in K_*^{\alpha\varepsilon,h_\varepsilon F}(\tilde A\rtimes G)$ such that $p_j(x)=x_j$ where $p_j$ is the composition of the projection on the $j^{th}$ factor $\tilde A \rtimes G \rightarrow K_A \rtimes G$ and the Morita equivalence in K-theory.
\[\begin{tikzcd}
RK^G(P_E(G),\tilde A) \arrow{r}{\mu_{G,\tilde A}^{\alpha\varepsilon,E, h_\varepsilon F}} \arrow{rd}{\mu_{G,\tilde A}^E} & K^{\alpha\varepsilon,h_\varepsilon F}(\tilde A\rtimes G) \arrow{d}{\iota_{\alpha\varepsilon, h_\varepsilon F}}\\
                                                      \                          & K(\tilde A\rtimes G)
\end{tikzcd}\]

If $\iota_{\alpha\varepsilon,h_\varepsilon F}(x)$ is in the range of $\mu_{G,\tilde A}$, there exists $E\in\mathcal E$ and $z\in RK^G(P_E(G),\tilde A)$ such that $\mu_{G,\tilde A}^E(z)=\iota_{\alpha\varepsilon,h_\varepsilon F}(x)$. Denote by $(z_j)$ the element of $\prod_j RK^G(P_E(G),A)$ corresponding to $z$ under the isomorphism $RK^G(P_E(G),\tilde A) \cong \prod_j RK^G(P_E(G),A)$ of lemma \ref{prod}. By remark \ref{approximation}, there exists $F''\in\mathcal E$ such that $h_\varepsilon F\subseteq F''$ and such that $\mu_{G,\tilde A}^{\lambda \alpha \varepsilon, E,F''}(z)=\iota_{\alpha\varepsilon,h_\varepsilon F}^{\lambda\alpha\varepsilon,F''}(x)$. \\

By naturality, $\mu_{G, A}^E(z_j)=0$. As $G$ admits a cocompact example for $\underline E G$, there exists $E'\in \mathcal E$ such that $E\subseteq E'$ and $q_E^{E'}(z)=0$. Since \[\mu_{G,\tilde A}^{\lambda\alpha\varepsilon, E,F''}(z)= \mu_{G,\tilde A}^{\lambda\alpha\varepsilon, E',F''}\circ q_E^{E'}(z),\]
we have that $\iota_{\alpha\varepsilon,h_\varepsilon F}^{\lambda\alpha\varepsilon,F''}(x)=0$ in $K^{\lambda\alpha\varepsilon,F''}(\tilde A\rtimes_r G)$. By naturality, for any $i$ such that $F''\subseteq F_i$, $\iota_{\varepsilon, F}^{\lambda\alpha\varepsilon,F_i}(x_i)=0$ is satisfied, which contradicts our assumption.\\  
%\iota_{\alpha\varepsilon,h_\varepsilon F}^{\lambda\varepsilon,F_j}(x_j)\neq 0$ and $\mu_{G,A}(z)=0$ so that $\mu_{G,A}$ is not injective. Hence either $\mu_{G,\tilde A}$ is not onto, either $\mu_{G,A}$ is not one to one. \\
\qed
\end{dem}

\begin{rk}
The theorems \ref{Quant1} and \ref{Quant2} provide examples that satisfy (PAP). Recall that every a-T-menable groupoid satisfies the Baum-Connes conjecture with coefficients, hence, if $G^{(0)}$ is compact and $G$ admits a cocompact example for $\underline E G$, $G$ satisfies the hypothesis of the theorem \ref{PAPG}.
\end{rk}

\begin{rk}
The theorem \ref{PAPG} provides an obstruction for the Baum-Connes conjecture to be satisfied. 
\end{rk}

%%%%%%%%%%%%%%%%%%%%%%%%%%%

%%%%%%%%%%%%%%%%%%%%%%
\section{Controlled assembly maps for coarse spaces}

In this section, $X$ will be a discrete metric space with bounded geometry, and $\mathcal E$ is the coarse structure generated by its controlled subsets. We also fix a separable Hilbert space $H$. For $R>0$, $\Delta_R$ is $\{(x,y)\in X\times X\text{ s.t. }d(x,y)<R\}$. We construct a controlled assembly map for coarse space $(X,\mathcal E)$ in the same way as for groupoids.  \\

%Recall first the construction of the Roe algebra of $X$ with coefficients in a $C^*$-algebra $B$, which can be found in \cite{SkTuYu}. 
$H_B$ denotes the standard $B$-Hilbert module $H\otimes B$. Recall that for every $x,y\in X$, and $T\in\mathcal L_B(l^2(X)\otimes H_B)$, we put $T_{xy}\in\mathcal L_B(H_B)$ to be the unique operator such that $\langle T_{xy}\xi,\eta\rangle = \langle T(e_x\otimes \xi),e_y\otimes\eta\rangle $ for every $x,y\in X$ and every $\xi,\eta\in H_B$.\\
%$T_{xy}= \chi_y T\chi_x$ , where $\chi_x,\chi_y$ are the characteristic functions of $\{x\}$ and $\{y\}$, seen as projection operators.\\

%For any positive number $R>0$, define the family of linear subspaces 
%\[C_R[X,B]=\{T\in \mathcal L(l^2(X)\otimes H_B) \text{ s.t. } T_{xy}\in \mathfrak K(H_B) \text{ and }T_{xy}=0 \text{ for }(x,y)\not\in \Delta_R  \}\]
%and $C^*(X,B)$ is the completion of $\cup_{R>0} C_R[X,B]$ for the operator norm in $\mathcal L(l^2(X)\otimes H_B) $. 

Remark that the $C^*$-algebra $C^*(X,B)$ is filtered by $\mathcal E$, and also by $\R_+^*$, seen as a coarse structure. Indeed, the composition law $R\circ R'= R+R'$ provides $\R_+^*$ with a coarse structure and $C_R[X,B] C_{R'} [X,B]\subseteq C_{R+R'}[X,B]$. For the $\mathcal E$-filtration, one can put :
\[C_E[X,B]=\{T\in \mathcal L(l^2(X)\otimes H_B) \text{ s.t. } T_{xy}\in \mathfrak K(H_B) \text{ and }T_{xy}=0 \text{ for }(x,y)\not\in E \}\quad \forall E\in\mathcal E.\]

To construct the coarse assembly map, we will need the following proposition.\\
 
Let $A$ be a $C^*$-algebra. The theorem \ref{Xfunctor} allows us to take the image of the exact sequence $0 \rightarrow SA \rightarrow CA \rightarrow A \rightarrow 0 $ under the functor $C^*(X,\cdot)$ to get the following filtered semi-split exact sequence 
\[0 \rightarrow C^*(X,SA) \rightarrow C^*(X,CA) \rightarrow C^*(X,A) \rightarrow 0.\] 
Let $D_{X,A} : \hat K_*(C^*(X,A))\rightarrow \hat K_*(C^*(X,SA))$ be the controlled boundary morphism associated to this last extension.

\begin{prop}\label{InverseEven}
Let $A$ be a $C^*$-algebra. 
%The theorem \ref{Xfunctor} allows us to take the image of the exact sequence $0 \rightarrow SA \rightarrow CA \rightarrow A \rightarrow 0 $ under the functor $C^*(X,\cdot)$ to get the following filtered semi-split exact sequence 
%\[0 \rightarrow C^*(X,SA) \rightarrow C^*(X,CA) \rightarrow C^*(X,A) \rightarrow 0.\] 
%Let $D_{X,A} : \hat K_*(C^*(X,A))\rightarrow \hat K_*(C^*(X,SA))$ the controlled boundary morphism associated to this last extension. 
Then there exists a control pair $(\lambda,h)$, independent of $X$ and $A$, such that $D_{X,A}$ is $(\lambda,h)$-invertible.
\end{prop}

\begin{dem}
Recall that $0\rightarrow \mathfrak K(l^2(\N)) \rightarrow \mathcal T_0\rightarrow S\rightarrow 0 $ is the Toeplitz extension. Let $\Psi$ be the obvious $*$-homomorphism $SC^*(X,A)\rightarrow C^*(X,SA) $. The following diagram has exact rows and commutes
\[\begin{tikzcd} 
SC^*(X,A) \arrow{d}{\Psi}\arrow{r} & CC^*(X,A) \arrow{d}\arrow{r} & C^*(X,A) \arrow{d} \\ 
C^*(X,SA) \arrow{r}          & C^*(X,CA) \arrow{r}          & C^*(X,A)
\end{tikzcd}\]
where vertical arrows are obvious inclusions. Remark \ref{rk3.8} implies that 
\[D_{X,A} = \Psi_*\circ D_{C^*(X,A)}.\] 
The following diagram also has exact rows and commutes
\[\begin{tikzcd} 
\mathfrak K(l^2(\N))\otimes C^*(X,A) \arrow{d}\arrow{r} & \mathcal T_0 \otimes C^*(X,A) \arrow{d}\arrow{r} & SC^*(X,A) \arrow{d}{\Psi} \\ 
C^*(X,\mathfrak K(l^2(\N))\otimes A) \arrow{r}          & C^*(X,\mathcal T_0 \otimes A) \arrow{r}          & C^*(X,SA)
\end{tikzcd}\]
where vertical arrows are obvious inclusions. Remark \ref{rk3.8} implies that 
\[D_{\mathfrak K(l^2(\N))\otimes C^*(X,A),\mathcal T_0\otimes C^*(X,A)} = D_{C^*(X,A),C^*(X,\mathcal T_0\otimes A)}\circ\Psi_*.\] 
A simple computation shows that 
\[D_{C^*(X,A),C^*(X,\mathcal T_0\otimes A)}\circ D_{X,A} \sim D_{\mathfrak K(l^2(\N))\otimes C^*(X,A),\mathcal T_0\otimes C^*(X,A)}\circ D_{C^*(X,A)} \sim \mathcal M_{C^*(X,A)}\]
\qed
\end{dem}

\begin{rk}\label{rkInverse} This result induces a similar statement in $K$-theory. %and in $KK$-theory similar statements. Namely : 
Namely, the boundary maps of the extensions $0 \rightarrow C^*(X,SA) \rightarrow C^*(X,CA) \rightarrow C^*(X,A) \rightarrow 0$ and
$0 \rightarrow\mathfrak K(l^2(\N))\otimes C^*(X,A) \rightarrow \mathcal T_0 \otimes C^*(X,A) \rightarrow SC^*(X,A) \rightarrow 0$ are inverse of each other in $K$-theory.
%\begin{itemize}
%\item[$\bullet$] the boundary maps of the extensions $0 \rightarrow C^*(X,SA) \rightarrow C^*(X,CA) \rightarrow C^*(X,A) \rightarrow 0$ and
%$0 \rightarrow\mathfrak K(l^2(\N))\otimes C^*(X,A) \rightarrow \mathcal T_0 \otimes C^*(X,A) \rightarrow SC^*(X,A) \rightarrow 0$ are inverse of each other in $K$-theory,  
%\item[$\bullet$] $[\partial_{K(l^2(\N))\otimes C^*(X,A), T_0 \otimes C^*(X,A)}]$ and $[\partial_{C^*(X,SA),C^*(X,CA)}]$ are $KK$-inverse of each other.
%\end{itemize}
\end{rk}

%%%%%%%%%%%%%%%%%%%%%%%%%%%%%%%%%%%%%%%%
\subsection{Controlled Roe transform}
%%%%%%%%%%%%%%%%%%%%%%%%%%%%%%%%%%%%%%%%

Every $K$-cycle $z\in KK(A,B)$ can be represented as a triplet $(H_B, \pi, T)$ where :
\begin{itemize}
\item[$\bullet$]$\pi : A\rightarrow \mathcal L_B(H_B)$ is a $*$-representation of $A$ on $H_B$.
\item[$\bullet$]$T\in \mathcal L_B(H_B)$ is a self-adjoint operator.
\item[$\bullet$] $T$ and $\pi$ satisfy the $K$-cycle condition, i.e. $[T,\pi(a)]$, $\pi(a)(T^*-T)$ and $\pi(a)(T^2-id_{H_B})$ are compact operators in $\mathfrak K_B(H_B)\cong \mathfrak K \otimes B$ for all $a\in A$.\\
\end{itemize}

We first define a controlled morphism $\hat \sigma_X(z) : \hat K(A)\rightarrow \hat K(B)$ of every $z\in KK(A,B)$, which we name the controlled Roe transform. It induces $-\otimes \sigma_X(z)$ in $K$-theory, and will be needed in the definition of the controlled coarse assembly map. Recall that if $\phi : A \rightarrow B$ is a $*$-homomorphism, we denote by $\phi_X : C^*(X,A)\rightarrow C^*(X,B)$ the induced $*$-homomorphism.

%\subsubsection{Odd case} %%%%%%%%%%%%%%

For $z\in KK_1(A,B)$, represented by $(H_B,\pi,T)\in E(A,B)$, define $P=(\frac{1+T}{2})\in \mathcal L_B(H_B)$ and 
%$P_X\in\mathcal L(H_{C^*(X,B)})$, and  
\[E^{(\pi,T)} = \{(a,P\pi(a)P + y) : a\in A,y\in  B\otimes \mathfrak K\} \]
which is a $C^*$-algebra such that the following sequence :
\[\begin{tikzcd}[column sep = small]0\arrow{r} & B\otimes \mathfrak K \arrow{r} & E^{(\pi,T)}\arrow{r} & A\arrow{r} & 0 \end{tikzcd}.\]
is exact and semi-split by the completely positive section $s : A\rightarrow B\otimes\mathfrak K ; a\mapsto P\pi(a)P$. Define $E_X = C^*(X,E^{(\pi,T)})$. Up to the $*$-isomomorphism $C^*(X,B\otimes\mathfrak K)\cong C^*(X,B)$, the following sequence
\[\begin{tikzcd}[column sep = small]0\arrow{r} & C^*(X,B) \arrow{r} & E_X^{(\pi,T)}\arrow{r} & C^*(X,A)\arrow{r} & 0 \end{tikzcd}.\]
is exact and semi-split by the completely positive section $s_X : C^*(X,A)\rightarrow E_X^{(\pi,T)}$.\\

The same proofs as Propositions \ref{ClassIndepedance}, \ref{Kasparov1} and \ref{Kasparov} yield the following results.

\begin{prop}
The controlled boundary map $D^{(\pi,T)}=D_{C^*(X,B),E_X^{(\pi,T)}}$ of the extension $E_X^{(\pi,T)}$ only depends on the class $z$.
\end{prop}

\begin{definition}
For every $z=[H_B,\pi,T]\in KK_1(A,B)$, we define the Roe transformation $\hat\sigma_X$ as 
\[\hat\sigma_X(z)= D_{C^*(X,B),E_X^{(\pi,T)}}\quad.\]
It is a $(\alpha_D,k_D)$-controlled morphism $\hat K(C^*(X,A))\rightarrow \hat K(C^*(X,B))$ of odd degree.\\

Let $z\in KK(A,B)$ be an even $K$-cycle. Recall that $[\partial_{SB}]\in KK_1(B,SB)$ is the $K$-cycle implementing the boundary of the extension $0\rightarrow SB\rightarrow CB\rightarrow B\rightarrow 0$, and $[\partial]\in KK_1(\C,S)$ is the Bott generator. Recall from proposition \ref{InverseEven} that $D_{X,A}$  and $D_{ C^*(X,A),C^*(X,\mathcal T_0\otimes A) }$ are controlled inverse of each other. We will denote $D_{ C^*(X,A),C^*(X,\mathcal T_0\otimes A) }$ by $T_{X,A}$.\\

As $z\otimes_B [\partial_{SB}]$ is an odd $K$-cycle, we define
\[\hat\sigma_X(z):= T_{X,B}\circ \hat\sigma_X(z\otimes[\partial_{SB}]).\] 
\end{definition}

The controlled Roe transform satisfies the following.

\begin{prop}\label{Roe2}
Let $A$ and $B$ two $C^*$-algebras. For every $z\in KK_*(A,B)$, there exists a control pair $(\alpha_X,k_X)$ and a $(\alpha_X,k_X)$-controlled morphism
\[\hat\sigma_X(z) : \hat K(C^*(X,A))\rightarrow \hat K(C^*(X,B))\]
of the same degree as $z$, such that
\begin{enumerate}
\item[(i)] $\hat\sigma_X(z)$ induces right multiplication by $\sigma_X(z)$ in $K$-theory ;
\item[(ii)] $\hat\sigma_X$ is additive, i.e.
\[\hat\sigma_X(z+z')=\hat\sigma_X(z)+\hat\sigma_X(z').\]
\item[(iii)] For every $*$-homomorphism $f : A_1\rightarrow A_2$,
\[\hat\sigma_X(f^*(z))=\hat\sigma_X(z)\circ f_{X,*}\] for all $z\in KK_*(A_2,B)$.
\item[(iv)] For every $*$-homomorphism $g : B_1\rightarrow B_2$,
\[\hat\sigma_X(g_*(z))= g_{X,*}\circ \hat\sigma_X(z)\] for all $z\in KK_*(A,B_1)$.
\item[(v)] $\hat\sigma_X([id_A]) \sim_{(\alpha_X,k_X)} id_{\hat K(C^*(X,A))}$,
\item[(vi)] Let $0\rightarrow J\rightarrow A\rightarrow A/J\rightarrow 0$ be a semi-split extension of $C^*$-algebras and $[\partial_J]\in KK_1(A/J,J)$ be its boundary element. Then 
\[\hat\sigma_X([\partial_{J,A}])=D_{C^*(X,J),C^*(X,A)}.\] 
\end{enumerate}
\end{prop}

We now show that the Roe transform respects in a quantitative way the Kasparov product. Let us recall the following result from \cite{lafforgue2002k}. It states that every $KK$-element comes from the product of an element coming from a $*$-homomorphism and an element coming from the inverse in $KK$-theory of a $*$-homomorphism. The following lemma is a particular case of decomposition property $d$, defined in \ref{DecompositionPropertyD}.

\begin{lem}[\cite{lafforgue2002k}, lemma $1.6.11$] Let $A$ and $B$ be two $C^*$-algebras and $z\in KK_0(A,B)$. Then, there exists a $C^*$-algebra $A_1$, an element $\alpha \in KK(A,A_1)$ and $*$-homomorphims $\theta : A_1 \rightarrow A$ and $\eta : A_1 \rightarrow B$ such that
$\theta^*(\alpha) = id_{A_1}$, $\theta_*(\alpha) = id_{A}$ and $\theta^*(z) = \eta$.
\end{lem}

\begin{prop} There exists a control pair $(\alpha_R,k_R)$ such that for every $C^*$-algebras $A$, $B$ and $C$, and every $z\in KK(A,B),z'\in KK(B,C)$, the controlled equality
\[\hat\sigma_X(z\otimes_B z') \sim_{\alpha_R,k_R} \hat\sigma_X(z')\circ \hat\sigma_X(z)\]
holds.
\end{prop}

\begin{dem}
Assume $\alpha\in KK_0(A,B)$. By naturality, the previous lemma reduces the proof to the special case of $\alpha$ being the inverse of a $*$-homomorphism $\theta : B\rightarrow A$ in $KK$-theory : $\alpha\otimes_B [\theta]=1_A$. Let $z\in KK(B,C)$ :
\[\begin{array}{rcl}
\hat\sigma_X (\alpha\otimes z) & \sim_{\alpha_J^2,k_J*k_J} &  \hat\sigma_X (\alpha\otimes z)\circ \hat\sigma_X (\alpha\otimes [\theta]) \\
			& \sim & \hat\sigma_X (\alpha\otimes z)\circ \hat\sigma_X (\theta_*(\alpha))\\
			& \sim & \hat\sigma_X (\alpha\otimes z)\circ \theta_{X,*}\circ \hat\sigma_X (\alpha)\\
			& \sim & \hat\sigma_X (\theta^*(\alpha\otimes z))\circ \hat\sigma_X (\alpha)\\
			& \sim & \hat\sigma_X (z)\circ \hat\sigma_X (\alpha) \\
\end{array}\] 
because $\theta^*(\alpha\otimes z)=\theta^*(\alpha)\otimes z=1\otimes z =z$. The control on the propagation of the first line follows from remark \ref{rk2.5} and point $(v)$, the other lines are equal by points $(iii)$ and $(iv)$, hence $(\alpha_R,k_R)$ can be taken to be $(2 \alpha_X^{4},( k_X)^{*2})$. If $z'$ is even, we can apply the same argument.\\

Let $z$ and $z'$ be odd $KK$-elements. Then :
\[\begin{array}{rcl}
\hat\sigma_X (z\otimes z') & = &  \hat\sigma_X (z\otimes_B [\partial_{B}]\otimes_{SB} [\partial_B]^{-1}\otimes_B z') \\
			& \sim & \hat\sigma_X ( [\partial_B]^{-1}\otimes_B z')\circ \hat\sigma_X (z\otimes_B [\partial_{B}])\\
			& \sim & \hat\sigma_X ( [\partial_B]^{-1}\otimes_B z')\circ \hat\sigma_X ( [\partial_B])\circ\hat\sigma_X( [\partial_B]^{-1})\circ \hat\sigma_X (z\otimes_B [\partial_{B}])\\		
			& \sim & \hat\sigma_X (  z')\circ \hat\sigma_X (z),\\
\end{array}\] 
where we used the previous case for the second line, Lemma \ref{Roe2} for the third line, and Proposition \ref{InverseEven} for the last one.\\
\qed
\end{dem}

\subsection{Controlled coarse assembly maps}

%If $(X,\mathcal E_X)$ is a coarse space, and $E\in\mathcal E_X$ a controlled subset, any simplex $\eta$ of the Rips complex $P_E(X) = \{m \in Prob(X)\text{ s.t. supp }m \subseteq E\}$ can be written as $\eta = \sum_{x\in X} \lambda_x(\eta) \delta_x$, where $\delta_x$ si the Dirac probability at $x$, and $\lambda_x : P_E(X)\rightarrow [0,1]$ is a continuous function. Set :
Let $E\in\mathcal E$ be a controlled subset. Then any probability $\eta$ of the Rips complex $P_E(X)$ can be written as $\eta = \sum_{x\in X} \lambda_x(\eta) \delta_x$, where $\delta_x$ si the Dirac probability at $x$, and $\lambda_x : P_E(X)\rightarrow [0,1]$ is a continuous function. Set :
\[ h_E : \left\{\begin{array}{rcl} X \times X & \rightarrow & C_0(P_E(X))\\  (x,y) & \mapsto & \lambda_x^{\frac{1}{2}}\lambda_y^{\frac{1}{2}}\end{array}\right. \]  
Let $(e_x)_{x\in X}$ be the canonical basis of $l^2(X)$, $e$ be a rank-one projection in $H$ and $P_E$ be defined as the extension by linearity and continuity of
\[P_E(e_x\otimes\xi\otimes f)= \sum_{y\in X} e_y\otimes (e\xi)\otimes (h(x,y)f)\] 
for every $x\in X$, $\xi\in H$ and $f\in C_0(P_E(X))$. As $\sum_{x\in X} \lambda_x =1$, $P_E$ is a projection of $\mathfrak K(l^2(X)) \otimes C_0(P_E(X))$ of controlled support : $\text{supp }P_E\subseteq E$. Indeed, $\lambda_x^{\frac{1}{2}}\lambda_y^{\frac{1}{2}} =0$ as soon as $(x,y)\notin E$. Hence $P_E$ defines a class $[P_E,0]_{\varepsilon, E'}\in K_0^{\varepsilon, E'} (C^*(X,C_0(P_E(X)))$ for any $\varepsilon\in (0,\frac{1}{4})$ and any $E'\in\mathcal E$ satisfying $E\subseteq E'$.\\

For every $C^*$-algebra $B$ and every controlled subsets $E,E'\in\mathcal E$ such that $E\subseteq E'$, the canonical inclusion $P_E(X)\hookrightarrow P_{E'}(X)$ induces a $*$-homomorphism $q_E^{E'} : C_0(P_{E'}(X))\rightarrow C_0(P_{E}(X))$, hence a map $(q_E^{E'})^* : KK(C_0(P_E(X)),B)\rightarrow KK(C_0(P_{E'}(X)),B)$ in $KK$-theory. It induces another map $((q_E^{E'})_X)_* : K(C^*(X,C_0(P_{E'}(X))))\rightarrow K(C^*(X,C_0(P_{E}(X))))$ in $K$-theory. The family of projections $P_E$ are compatible with the morphisms $q_E^{E'}$, i.e. $((q_E^{E'})_X)_*[P_{E'},0]_{\varepsilon,E'} = [P_{E},0]_{\varepsilon,E}$, for every $\varepsilon\in (0,\frac{1}{4})$.
%Moreover, the inclusion being an isometry, we have a map $K(C^*(P_E(X), B)) \rightarrow K(C^*(P_E(X), B))$, still denoted $q_E^{E'}$.

\begin{definition}
Let $B$ a $C^*$-algebra, $\varepsilon\in (0,\frac{1}{4})$ and $E,F\in\mathcal E_X$ controlled subsets such that $k_X(\varepsilon).E\subseteq F$. The controlled coarse assembly map $\hat\mu_{X,B}=(\mu_{X,B}^{\varepsilon,E,F})_{\varepsilon,E}$ is defined as the family of maps
\[\hat\mu_{X,B}^{\varepsilon, E,F} :\left\{\begin{array}{rcl} KK(C_0(P_E(X)),B) & \rightarrow & K^{\varepsilon, F}(C^*(X,B)) \\
					z & \mapsto & \iota_{\alpha_X \varepsilon',k_X(\varepsilon').F'}^{\varepsilon,F}\circ\hat\sigma_X(z)[P_{E},0]_{\varepsilon', F'}\end{array}\right.\]
where $\varepsilon'$ and $F'$ satisfy :
\begin{itemize}
\item[$\bullet$] $\varepsilon'\in (0,\frac{1}{4})$ such that $\alpha_X \varepsilon'\leq \varepsilon$,
\item[$\bullet$] and $F'\in\mathcal E$ such that $E\subseteq F'$ and $k_X(\varepsilon').F'\subseteq F$.
\end{itemize}
%are chosen not to exceed $\varepsilon$ and $E$ when composed with the propagation of the controlled morphisms. 
\end{definition}

\begin{rk} The controlled coarse assembly map is compatible with the structure morphisms $q_E^{E'}$. Indeed, for every $E,E'\in \mathcal E$ such that $E\subseteq E'$, by proposition \ref{Roe2}, 
\[\hat\sigma_X((q_E^{E'})^*(z))[P_{E'},0]_{\varepsilon,E'}  = \hat\sigma_X(z)\circ ((q_E^{E'})_X)_*[P_E',0]_{\varepsilon,E'}= \hat\sigma_X(z)[P_E,0]_{\varepsilon,E}.\] 
Hence $\hat\mu_{X,B}^{\varepsilon,E,F}\circ(q_E^{E'})^* =\hat\mu_{X,B}^{\varepsilon,E',F}$.
\end{rk}

\begin{rk} The controlled coarse assembly map is also compatible with the structure morphisms $\iota_{\varepsilon,E}^{\varepsilon',E'}$, i.e. $\iota_{\varepsilon,F}^{\varepsilon',F'}\circ\hat\mu_{X,B}^{\varepsilon,E,F} =\hat\mu_{X,B}^{\varepsilon',E,F'}$ for every $F\subseteq F'$ and $\varepsilon\leq \varepsilon'$ such that this equality is defined. 
\end{rk}

\begin{rk}According to Proposition \ref{Roe2}, $\hat\sigma_X(z)$ induces right-multiplication by $\sigma_X(z)$. Hence, the controlled coarse assembly map $\hat \mu_{X,B}$ induces the coarse assembly map $\mu_{X,B}$ in $K$-theory.
\end{rk}

\begin{rk}
This assembly map is defined for the usual Roe algebra of $X$, but could be defined for any "nice" completion of the algebraic Roe algebra $\cup_{E\in \mathcal E_X} C_E[X]$. In particular, we can define an assembly map with values in the controlled $K$-theory of the maximal Roe algebra $C_{max}
^*(X)$, that we will denote by $\hat \mu_{X}^{max}$.\end{rk}

%%%%%%%%%%%%%%%%%%%%%%

%%%%%%%%%%%%%%%%%%%%%%%%%
\section{Applications to Coarse Geometry}

We present in this section a result on the equivalence between the controlled assembly map for a discrete metric space with bounded geometry $X$ with coefficients in a $C^*$-algebra $B$ and the controlled assembly map for the coarse groupoid $G(X)$ with coefficients in the $G(X)$-algebra $l^\infty(X,B\otimes \mathfrak K)$. This result is applied to show that any such space that admits a fibred coarse embedding into Hilbert space satisfies the maximal controlled Baum-Connes conjecture.

\subsection{Equivalence between the controlled coarse assembly map for $X$ and the controlled assembly map for $G$ with coefficients in $l^\infty(X,\mathfrak K)$}

In this section, we prove how the result of G. Skandalis, J.-L. Tu and G. Yu \cite{SkTuYu} extends to the setting of controlled $K$-theory. \\

Recall from Lemma $4.4$ in \cite{SkTuYu} that, for every $C^*$-algebra $B$, there exists a natural isomorphism of $C^*$-algebras 
\[\Psi_B : l^\infty(X,B\otimes\mathfrak K)\rtimes_r G(X)\rightarrow C^*(X,B).\]
Moreover, it is filtered in the strong sense : for every entourage $E\subseteq X\times X$, $\Psi_B(C_{\overline E}(G,B))= C_E[X,B]$.\\

The following theorem is proved in \cite{SkTuYu}. It states the equivalence between the coarse Baum-Connes conjecture with coefficients in $B$ and the Baum-Connes conjecture for $G(X)$ with coefficients in $l^\infty(X,B\otimes\mathfrak K)$.  

\begin{thm}[\cite{SkTuYu}]
Let $X$ be a discrete metric space with bounded geometry. Let $\Psi_B$ be the previous isomorphism, $x\in X$ and $\iota :\{x\}\rightarrow G(X)$ be the natural inclusion of groupoids. Denote by $G=G(X)$ the coarse groupoid of $X$ and by $\tilde B$ the $G$-algebra $l^\infty (X,B\otimes\mathfrak K)$. Then, for every controlled subset $E\subseteq X\times X$, the following diagram is commutative with vertical arrows being isomorphisms :
\[\begin{tikzcd}
RK_*^G(P_{\overline{E}}(G),\tilde B) \arrow{r}{\mu_{G,\tilde B}^{\overline E}}\arrow{d}{\iota^*}& K_*(\tilde B\rtimes_r G)\arrow{d}{(\Psi_B)_*}\\
RK_*(P_E(X),B) \arrow{r}{\mu_{X,B}^E}& K_*(C^*(X,B))
\end{tikzcd},\]
where $\iota^*$ is the natural transformation induced by $\iota$ and $d= \sup_E d$.
\end{thm}

We shall prove a controlled analogue of this result which induces it in $K$-theory. We need the following lemmas. %s.

%\begin{lem}
%Let $G$ be an étale groupoid,$x\in G$, $Z$ a proper $G$-space and $B$ a $C^*$-algebra. Denote by $\tilde A$ the $G$-algebra $C_0(Z)$ and $\iota : \{x\} \rightarrow G$ the natural inclusion of groupoids. Then :
%\[\iota^* : RK^G(Z,l^\infty_B)\rightarrow KK(A_x,B)\]
%is an isomorphism of $\Z_2$-graded abelian groups. 
%\end{lem}

%\begin{dem}
%We define an inverse for $\iota^*$ : for $z=[H_{B_x},\pi,T]\in KK(A_x,B_x)$, define $\eta(z)= [H_B,\tilde\pi,\tilde T]$ where 
%\[(\tilde\pi) = \pi \otimes id\]
%\qed
%\end{dem}
\begin{lem}[Lemma $4.7$ \cite{SkTuYu}]\label{iota}
Let $x\in X$ and $\iota : \{x\}\hookrightarrow G$ the natural inclusion of groupoids. Then 
\[\iota^* : KK^G(C_0(P_{\overline E}(G),\tilde B) \rightarrow KK(C_0(P_{E}(X),B) \]
is an isomorphism of $\Z_2$-graded abelian groups.
\end{lem}

The reader can find a proof in \cite{SkTuYu} (Lemma $4.7$). We recall the explicit construction of the inverse 
\[j:KK(C_0(P_{E}(X),B) \rightarrow KK^G(C_0(P_{\overline E}(G),\tilde B)\] 
of $\iota^*$. Let $(H_B,\pi,T)\in\mathbb E (C_0(P_{E}(X),B)$ be a standard $K$-cycle. Let $\tilde B= l^\infty(X,B\otimes\mathfrak K)$ seen as Hilbert module over itself, $(\tilde \pi (a)\xi)(x) = \pi(a(x))\xi(x)$ and $(\tilde T\xi ) (x) = T\xi(x)$, for every $x\in X$ and $\xi\in E$. Then $j([H_B,\pi,T])=[\tilde B,\tilde \pi, \tilde T]$.

\begin{lem} Let $E\subseteq X\times X$ be controlled subset and $B$ be $C^*$-algebra. Denote $C_0(P_E(X))$ by $A$, $C_0(P_{\overline E}(G))$ by $\tilde A$ and $\tilde B = l^\infty(X,B\otimes \mathfrak K)$. Then, for every $z\in KK^G(\tilde A,\tilde B)$, the following equality of controlled morphisms holds :
\[\hat\sigma_X(\iota^*(z))\circ (\Psi_A)_* = (\Psi_B)_*\circ \hat J_G(z).\]  
\end{lem}

%%% NEW NEW PROOF
\begin{dem}
Let $z\in KK_1^G(\tilde A,\tilde B)$. Let $\iota^*(z)$ be represented by the $K$-cycle $[H_{B},\pi,T]\in\mathbb E(A,B)$, and let $P=\frac{1+T}{2}$. Denote by $(\tilde B,\tilde \pi ,\tilde T)\in \mathbb E^G(\tilde A,\tilde B)$ the representative of $j(\iota^*(z))=z$ constructed as in lemma \ref{iota} and $\tilde P=\frac{1+\tilde T}{2}$. Recall that 
\[E^{(\pi,T)} = \{(x,P\pi(x)P+y : x\in A,y\in B\otimes\mathfrak K\},\] 
and $E^{(\pi,T)}_X=C^*(X,E^{(\pi,T)})$. \\

First, notice that $z$ is the boundary element in $KK^G(\tilde A,\tilde B)$ of the following extension 
\[0 \rightarrow \tilde B \rightarrow E'\rightarrow \tilde A \rightarrow 0\]
where $E'$ is the $G$-algebra $\{ (a,\tilde P \tilde \pi(a) \tilde P+y) : a\in \tilde A, y \in \tilde B  \}\subseteq \tilde A\oplus \mathcal M(\tilde B) $, and the $*$-homomorphisms are the obvious ones. Set 
\[E'_G=\{ (a,\tilde P_G \tilde \pi_G(a) \tilde P_G+y) : a\in \tilde A\rtimes_r G, y \in \tilde B\rtimes_r G    )\}.\] 
We take the previous extension under the reduced crossed product to get the following extension
\[0 \rightarrow \tilde B\rtimes_r G \rightarrow E_G'\rightarrow \tilde A\rtimes_r G \rightarrow 0.\]
By \ref{Kasparov1}, $J_G(z)$ is given by the controlled boundary of $E'_G$. \\

We shall define a $*$-homomorphism from $E'\rtimes G$ to $E^{(\pi,T)}_X$ that intertwines the two extensions. Extend the $*$-isomorphism $\Psi_B : \tilde B \rtimes_r G \rightarrow C^*(X,B)$ to $\tilde \Psi_B : \mathcal M(\tilde B \rtimes_r G) \rightarrow \mathcal M(C^*(X,B))$. Set $\Psi_{E'} (a,y) = (\Psi_A(a),\tilde \Psi_B(y)) $ for every $(a,y)\in E'_G$. This map makes the following diagram commutes
\[
\begin{tikzcd}[column sep = small]
0\arrow{r} & \tilde B\rtimes_r G \arrow{r} \arrow{d}{\Psi_B} & E'_G \arrow{r}\arrow{d}{\Psi_{E'}} &
			 \tilde A\rtimes_r G\arrow{r}\arrow{d}{\Psi_A} & 0 \\
0\arrow{r} & C^*(X,B) \arrow{r} & E^{(\pi,T)}_X  \arrow{r} & C^*(X,A)\arrow{r} & 0 
\end{tikzcd}.
\]
By remark \ref{rk3.8}, we get 
\[  (\Psi_B)_* \circ D_{\tilde B\rtimes_r G, E'_G} = D_{ C^*(X,B), E^{(\pi,T)}_X} \circ (\Psi_A)_*,\]
hence,
\[ (\Psi_B)_*\circ \hat J_G(z) = \hat\sigma_X(\iota^*(z)) \circ (\Psi_A)_*.\]
\qed
\end{dem}

%%% END NEWPROOF

\begin{thm}\label{BCCeq}
Let $B$ be a $C^*$-algebra, $E\in \mathcal E_X$ an entourage and $\overline E \in \mathcal E_G$ the corresponding compact open subset of $G$. With the above notations, for all $z\in RK^G(P_{\overline E}(G),\tilde B)$ and all $\varepsilon\in(0,\frac{1}{4})$, the following equality holds :
\[(\Psi_B)_*\circ\mu^{\epsilon,\overline E}_{G,\tilde B} (z) = \mu_{X,B}^{\epsilon,E}(\iota^*(z)).\]
\end{thm}

\begin{dem}
%Let $E$ be a compact subset of $G$ such that $\overline \Delta_R \subseteq E$.
By the previous lemma, we only need to check that $(\Psi_A)_*[\mathcal L_{\overline E},0]_{\varepsilon,\overline  E} = [P_{E},0]_{\varepsilon, E} $, which is trivial.\\
\qed
\end{dem}

\begin{rk}
This theorem remains true for the maximal version of the assembly map when $B=\C$. One then has to replace $\hat\mu_{G,\tilde \C}$ and $\hat\mu_{X}$ by $\hat\mu^{max}_{G,\tilde \C}$ and $\hat\mu^{max}_{X}$ respectively.
\end{rk}

This result induces the result of \cite{SkTuYu} in $K$-theory. It also implies interesting consequences for Coarse Geometry. Recall that if the groupoid $G$ satisfies the Baum-Connes conjecture with coefficients, it satisfies the controlled Baum-Connes conjecture. Interesting examples follow from the result of J-L. Tu \cite{TuThese} that a-$T$-menable groupoids satisfy the Baum-Connes conjecture with coefficients. In particular, \\

\begin{itemize}
\item[$\bullet$] amenable groupoids are a-$T$-menable.\\
\item[$\bullet$] Let $X$ be a uniformly discrete metric space with bounded geometry. Then, if $X$ is coarsely embeddable into a separable Hilbert space, $G(X)$ is a-$T$-menable \cite{SkTuYu}. \\
\end{itemize}

\subsection{Fibred coarse embedding}

We now present an application to fibred coarse embedding.

\begin{definition}
Let $X$ be a discrete metric space with bounded geometry and $B$ a $C^*$-algebra. We introduce the following properties.\\
\begin{itemize} 
\item[$\bullet$] $QI_{X,B}(E,E',F,\varepsilon)$ : for any $x\in KK(C_0(P_E(X)), B )$, then $\mu^{\varepsilon,E,F}_{X,B}(x) = 0$ implies $q_E^{E'}(x)=0$ in $KK^G(C_0(P_{E'}(X)),B)$.
\item[$\bullet$] $QS_{X,B}(E,F,F',\varepsilon,\varepsilon')$ : for any $y\in K^{\varepsilon,F}(C^*(X,B))$, there exists $x\in KK(C_0(P_E(G)),B)$ such that $\mu^{\varepsilon',E,F'}_{X,B}(x)=\iota_{\varepsilon,F}^{\varepsilon',F'}(y)$.\\
\end{itemize} 
Let $\lambda \geq 1$ be a positive number. We say that $X$ satisfies the controlled Baum-Connes conjecture with coefficients in $B$ with rescaling $\lambda$ if :
\begin{itemize} 
\item[$\bullet$] for every $\varepsilon \in (0,\frac{1}{4\lambda})$, every $E,F\in\mathcal E$ such that $k_X(\varepsilon).E\subseteq F$, there exists $E'\in \mathcal E$ such that $E \subseteq E'$ and $ QI_{X,B}(E,E',F,\varepsilon)$ holds; 
\item[$\bullet$] for every $\varepsilon \in (0,\frac{1}{4\lambda})$, every $F\in\mathcal E$, there exists $E,F'\in\mathcal E$ such that $k_X(\varepsilon).E \subseteq F'$ and $F\subseteq F'$ and $QS_{X,B}(E,F,F',\varepsilon,\lambda\varepsilon)$ holds. 
\end{itemize} 
If $\hat\mu_{X}$, is replaced by $\hat\mu^{max}_{X}$, we will say that $X$ satisfies the maximal controlled Baum-Connes conjecture with rescaling $\lambda$.\\
\end{definition}

Recall from Theorem $1$ in \cite{FinnSellFibred} that if $X$ admits a fibred coarse embedding into Hilbert space, then $G(X)_{|\partial \beta X}$ is a-T-menable. For interesting examples of this type, recall the definition of a box space. Let $\Gamma$ be a finitely generated group, and $\mathcal N$ a family of nested normal subgroups with trivial intersection, which have finite index in $\Gamma$. Take the coarse union of the quotients to construct a coarse space $X_{\mathcal N}(\Gamma)= \cup_{H\in \mathcal N } \Gamma/ H$. Then, $X_{\mathcal N}(\Gamma)$ admits a fibred coarse embedding if and only if $\Gamma$ is a-$T$-menable. But if $X_{\mathcal N}$ is an expander, it cannot be coarsely embedded into a Hilbert space, so just take an a-$T$-menable group which has a box space $X$ which is an expander to get a coarse space that is not coarsely embeddable into Hilbert space ($SL(2,\Z)$ for instance), but admits a fibred coarse embedding.\\

The last example gives the following corollary.

\begin{cor}\label{fibred}
Let $X$ be a coarse space that admits a fibred coarse embedding into Hilbert space. Then $X$ satisfies the maximal controlled Coarse Baum-Connes conjecture. %$\hat \mu_{X}^{max}$ is a controlled isomorphism, i.e. $X$ satisfies the controlled Coarse Baum-Connes conjecture.
\end{cor}

\begin{dem}
By theorem \ref{BCCeq}, it is sufficient to show that $G(X)$ satisfies the maximal controlled Baum-Connes conjecture with coefficients in $l^\infty(X,\mathfrak K)$. We will denote $l^\infty(X,\mathfrak K)$ by $l^\infty$.\\
  
The maximal crossed product turns restriction of a groupoid to invariant open subsets into exact sequences of $C^*$-algebras, hence  
\[0\rightarrow l^\infty \rtimes_{max} G_{|U} \rightarrow l^\infty \rtimes_{max} G \rightarrow l^\infty \rtimes_{max} G_{|Y} \rightarrow 0\]
is an exact sequence, with $Y=\partial\beta X$ and $U= Y^c$. Moreover 
\[l^\infty \rtimes_{max} G_{|U}\cong l^\infty_{|U} \rtimes_{max} G
\quad \text{ and } \quad l^\infty \rtimes_{max} G_{|Y}\cong (l^\infty/l^\infty_{|U}) \rtimes_{max} G,\] 
hence $[\partial_{l^\infty \rtimes G_{|U},l^\infty\rtimes_r G}]=j_G([\partial_{l^\infty_{|U},l^\infty}]) $. 
Recall from Proposition \ref{Kasparov1} that $J_G([\partial_{l^\infty_{|U},l^\infty}])=D_{l^\infty_{|U}\rtimes_r G,l^\infty \rtimes_r G}$, hence there exists a control pair $(\alpha,k)$ such that for every $z\in RK^G(P_E(G),l^\infty / l^\infty_{|U} )$, 
\[\mu_{G}^{\varepsilon,E,F}(z\otimes [\partial_{l^\infty_{|U},l^\infty}] ) 
\sim_{\alpha,k} D_{l^\infty_{|U}\rtimes_r G,l^\infty \rtimes_r G} \circ \mu_{G}^{\varepsilon,E,F}(z )\]

Hence the following diagram commutes :
\[\begin{tikzcd}
RK^G(P_E(G),l^\infty_{|Y}) \arrow{d}{\otimes[\partial_{l^\infty_{|U},l^\infty}]} \arrow{r}{\mu_{G}^{\varepsilon,E,F}} 
			& K_*^{\varepsilon,F}(l^\infty_{|Y} \rtimes_{max} G) \arrow{d}{D_{l^\infty_{|U}\rtimes_r G,l^\infty \rtimes_r G}} \\
RK^G(P_E(G),l^\infty_{|U} )\arrow{d}\arrow{r}{\mu_{G}^{\alpha\varepsilon,E,k(\varepsilon).F}} 
			& K_*^{\alpha\varepsilon,k(\varepsilon).F}(l^\infty_{|U} \rtimes_{max} G) \arrow{d} \\
RK^G(P_E(G),l^\infty)      \arrow{d}\arrow{r}{\mu_{G}^{\alpha\varepsilon,E,k(\varepsilon).F}} 
			& K_*^{\alpha\varepsilon,k(\varepsilon).F}(l^\infty \rtimes_{max} G)      \arrow{d} \\
RK^G(P_E(G),l^\infty_{|Y}) \arrow{d}{\otimes[\partial_{l^\infty_{|U},l^\infty}]}\arrow{r}{\mu_{G}^{\alpha\varepsilon,E,k(\varepsilon).F}} 
			& K_*^{\alpha\varepsilon,k(\varepsilon).F}(l^\infty_{|Y} \rtimes_{max} G) \arrow{d}{D_{l^\infty_{|U}\rtimes_r G,l^\infty \rtimes_r G}} \\
RK^G(P_E(G),l^\infty)               \arrow{r}{\mu_{G}^{\alpha\varepsilon,E,k(\varepsilon).F}} 
			& K_*^{\alpha\varepsilon,k(\varepsilon).F}(l^\infty_{|U} \rtimes_{max} G) \\
\end{tikzcd}.\]
Now, $G_{|Y}$ being a-T-menable and $G_{|U}$ being proper, $\mu_{G_{|Y},B}$ and $\mu_{G_{|U},B}$ are isomorphisms for any $G$-algebra $B$. By theorems \ref{Quant1} and \ref{Quant2}, the families of the four exterior horizontal maps satisfies the controlled Baum-Connes conjecture, and the controlled version of the five lemma concludes the proof.\\
\qed
\end{dem}

%%%%%%%%%%%%%%%%%%%%%%%%%%%%

\bibliographystyle{plain}
\bibliography{biblio2} 

\end{document}